\documentclass[preprint,12pt]{elsarticle}
\pdfoutput=1

\usepackage{graphicx}
\usepackage{amssymb}
\usepackage{amsthm}
\usepackage{graphicx}
\usepackage{graphics}
\usepackage{upgreek}
\usepackage{fixltx2e}
\usepackage{float}
\usepackage{epstopdf}
\usepackage{amsmath}
\usepackage{wrapfig}
\usepackage{pdflscape}
\usepackage{listings}
\usepackage{subcaption}
\usepackage{multirow}
\usepackage[percent]{overpic}
\usepackage{multicol}
\usepackage{color,soul}
\usepackage{stmaryrd}
\usepackage[hidelinks]{hyperref}
\usepackage{algorithmic}

\usepackage{booktabs}
\usepackage{epstopdf}
\usepackage[section]{placeins} 

\usepackage{calc}
\newlength\myheight
\newlength\mydepth
\settototalheight\myheight{Xygp}
\newcommand*\inlinegraphics[1]{%
  \settototalheight\myheight{Xygp}%
  \settodepth\mydepth{Xygp}%
  \raisebox{-\mydepth}{\includegraphics[height=\myheight]{#1}}%
}
\newcommand\orcid[1]{\href{https://orcid.org/#1}{\inlinegraphics{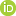}}}

\usepackage[titletoc]{appendix}

\usepackage[section]{algorithm}

\makeatletter
\def\BState{\State\hskip-\ALG@thistlm}
\makeatother

\journal{}

\graphicspath{ {Figs/} }

\newcommand{\etal}{et~al.}
\newcommand\floor[1]{\lfloor#1\rfloor}

\newcommand{\hfd}{\hat{f}^{\delta}}

\newcommand{\hf}{\hat{f}}
\newcommand{\hfI}{\hat{f}^{\delta I}}
\newcommand{\hfC}{\hat{f}^{\delta C}}
\newcommand{\hfD}{\hat{f}^{\delta D}}

\newcommand{\hud}{\hat{u}^{\delta}}
\newcommand\hudi[1]{\hat{u}^{\delta(#1)}}

\newcommand{\hu}{\hat{u}}
\newcommand{\huI}{\hat{u}^{\delta I}}

\newcommand\tb[1]{\tilde{\mathbf{#1}}}
\newcommand\hb[1]{\hat{\mathbf{#1}}}

\newcommand\px[2]{\frac{\partial #1}{\partial {#2}}}
\newcommand\pxi[3]{\frac{\partial^{#1}#2}{\partial {#3}^{#1}}}
\newcommand\dx[2]{\frac{\mathrm{d} #1}{\mathrm{d} #2}}
\newcommand\dxi[3]{\frac{\mathrm{d}^{#1}#2}{\mathrm{d} {#3}^{#1}}}
\newcommand\rint[2]{\int^{1}_{-1}{#1} \mathrm{d}{#2}}

\newcommand{\half}{\frac{1}{2}}

\newcommand\CONDITION[2]%

\begin{document}

\begin{frontmatter}

\title{Generalised Sobolev Stable Flux Reconstruction}

\author{W. Trojak\corref{cor1}~\orcid{0000-0002-4407-8956}}
\ead{wt247@cam.ac.uk}
\cortext[cor1]{Corresponding author}
\address{Department of Engineering, University of Cambridge, CB2 1PZ}

\begin{abstract}
	A new set of symmetric correction functions is presented for high-order flux reconstruction, that expands upon, while incorporating, all previous correction function sets and opens the possibility for improved performance. By considering FR applied to the linear advection equation, and through modification to the Sobolev norm, criteria are presented for a wider set of correction functions. Legendre polynomials are then used to fulfil these criterion and realise functions for third to fifth order FR. The sufficient conditions for the existence of the modified norms are also explored, before Fourier and Von Neumann analysis are applied to analytically find temporal stability limits for various Runge-Kutta temporal integration schemes. For all cases, correction functions are found that extend the temporal stability of FR. Two application-inspired investigations are performed that aim to explore the effect of aliasing and non-linear equations. In both cases unique correction functions could be found that give good performance, compared to previous FR schemes, while also improving upon the temporal stability limit.
	
\end{abstract}


\end{frontmatter}



\section{Introduction}
\label{sec:intro}

	High-order methods have steadily developed over the last half century based on the promise of allowing for complex problems to be tackled with potentially higher accuracy and greater efficiency. When the application of such methods to Large-Eddy Simulation (LES) is considered, as investigated by Chow~\cite{Chow2003} and Moin and Ghosal~\cite{Ghosal1996}, high-order becomes key in reducing grid requirements. This is due to reduction the global solution error that is seen in many cases when higher orders are used. But also high-order results in the decoupling of truncation error and aliasing error, which can allow for more effective use of sub-grid scale models.
	
	The early higher-order methods were mostly based on finite difference approaches, for example the fourth order finite difference method of Roberts and Weiss~\cite{Roberts1966}. Another such early high-order method introduced by Reed and Hill~\cite{Reed1973} was the Discontinuous Galerkin (DG) method. Through DG's own evolution --- see Cockburn~\etal~\cite{Cockburn2000} for a comprehensive review ---  other methods developed, including the Spectral collocation method~\cite{Quarteroni1990} and the Staggered Grid~(SG) approaches of Kopriva~\&~Kolias~\cite{Kopriva1996a} and Kopriva~\cite{Kopriva1996}. Huynh~\cite{Huynh2007}, and later Wang and Goa~\cite{Wang2009}, introduced a finite element method that built on the developments of DG and SG to form the foundation of what is now Flux Reconstruction (FR). FR is a high-order unstructured finite element type approach that uses the propagation of corrections between elements to form a piecewise continuous flux. The advantage of FR, compared to some other high-order approaches, is its suitability to current trends in heterogeneous computing --- in particular to many core architectures --- which is largely due to its use of domain subdivision~\cite{Vincent2017}.  
	
	The character of FR, including the wave propagating performance and temporal integration stability, have been found to be somewhat dependent on the manner in which the correction is applied~\cite{Huynh2007}. Several methods emerged for applying corrections, initially via a single parameter family~\cite{Vincent2010} that incorporated the work of Huynh~\cite{Huynh2007} as well as nodal DG~\cite{Hesthaven2007} and the variation of the spectral difference method~\cite{Kopriva1996,Liu2006,Sun2007} that was set in FR by Jameson~\cite{Jameson2010}. This single parameter family of correction functions was successfully extended for advection-diffusion~\cite{Huynh2009,Castonguay2014} as well for several multi-dimensional element types~\cite{Huynh2011,Williams2014,Sheshadri2016}. We will henceforth call this single parameter set of FR correction functions 'Original Stable FR'~(OSFR). Through the continued theoretical developments of FR, further links have emerged between FR and linear filtered DG~\cite{Allaneau2011,DeGrazia2014,Zwanenburg2016,Mengaldo2016}. 
	
	More recently, theory drawn from finite and spectral element methods has been applied to extend the range of stable FR correction functions. We will refer to these schemes as Extended range Stable FR (ESFR)~\cite{Vincent2015}. Further analysis of these schemes has shown that ESFR has the potential to be more stable when confronted by aliasing errors~\cite{Vermeire2016a}, as well as having the potential to reduce oscillations caused by dispersion when FR is applied as implicit LES~\cite{Vermeire2017a}.
	
	In this paper, a further extension to the correction function set is presented that aims to incorporate the methods previously presented, as well as furthering the potential for improved temporal stability and improved dispersion and dissipation characteristics. The ultimate aim is that this approach will broaden the future ways in which scheme stability is sought, while providing potential methods for improved FR calculation performance through an expanded correction function definition.
	
	Before presenting what we will refer to as Generalised Sobolev Stable FR (GSFR), the current state of the art is established, this will include the description of OSFR and ESFR in Section~\ref{sec:fr}. In Section~\ref{sec:gsobolev} the methodology for deriving GSFR will then be laid out together with the criteria defining the correction functions. Section~\ref{sec:stability} will explore the limitations on stability, with specific evaluation of the GSFR criteria for several orders presented in Section~\ref{sec:GSFR}. Also in this section the uniqueness of the new functions relative to those of OSFR and ESFR will be shown. Theoretical assessment of the performance of newly presented set will be performed via Fourier/Von Neumann analysis in Section~\ref{sec:vn}. Lastly, in Section~\ref{sec:numerical}, a linear heterogeneous equation and Euler's equations will be solved using the GSFR method in order to determine if, when applied to equation sets of significance to engineering, GSFR can give reduced aliasing error and dissipation. Conclusion will be made in Section~\ref{sec:conclusions}.
\section{Flux Reconstruction Method}
\label{sec:fr}

\begin{figure}[tbhp]
		\centering
		\includegraphics[width=0.4\linewidth]{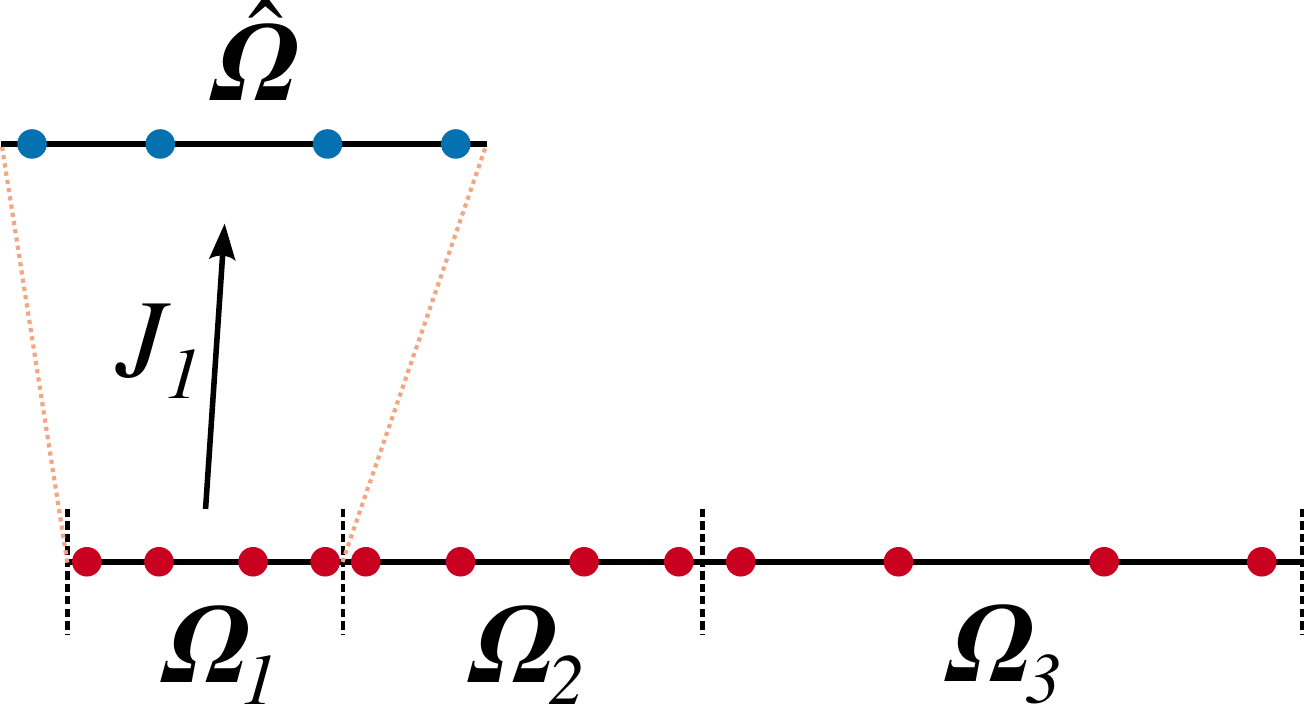}
		\caption{Diagram showing the mapping from the physical domain to the reference domain, with points to show the spatial discretisations.}
		\label{fig:FR_domain}
	\end{figure}

To introduce the numerical underpinning of FR, we begin by considering the 1D conservation law:
	\begin{equation}
		\px{u}{t} + \px{f}{x} = 0 \label{eq:linearConLaw}
	\end{equation}	
	where $u$ is the conserved quantity and $f=f(u)$ is some arbitrary flux function. The spatial domain of the solution will be taken as being $\mathbf{\Omega}$, with spatial variable $x\in\mathbf{\Omega}$. Flux reconstruction decomposes this domain into sub-domains, $\mathbf{\Omega}_i$, such that:
	\begin{equation}
		\mathbf{\Omega} = \bigcup_{j=1}^{m} \mathbf{\Omega}_j \quad \mathrm{and} \quad \bigcap^m_{i=0} \mathbf{\Omega}_j = \emptyset
	\end{equation}
	where $m$ is the number of sub-domains. Each sub-domain may then be transformed into a regularized reference domain, $\hb{\Omega}$, that is typically taken as being the closed interval $[-1,1]$ in 1D with reference spatial variable $\xi\in\hb{\Omega}$. These definitions now allow for the domain to be discretised such that in each sub-domain, $\mathbf{\Omega}_j$, there is an ordered set of points $\{x_0,\dots,x_{n_s}\}$. Within the reference domain, the discretisation points are the ordered set $\{\xi_0,\dots,\xi_{n_s}\}$. It is typical within FR to use a maximal order basis, which implies $n_s=(p+1)^d$. Here $p$ is the polynomial order, to be defined shortly, and $d$ is the dimensionality. The mapping from sub-domain to reference domain can then be characterised by the Jacobian, $J_j$, such that the transformed variables are:
	\begin{align}
		\hu_j &= J_ju_j\\
		\hf_j &= J_jf_j
	\end{align}
	Here superscript hats indicate that a variable has been transformed from the physical domain to the reference domain.	This set-up is further demonstrated in Fig.~\ref{fig:FR_domain}.
	
	Within each element a local polynomial fit of both the conserved variable and flux variable can then be defined by using the point values such that:
	\begin{align}
		\hud &= \sum_{i=0}^p \hat{u}(\xi_i)l_i(\xi) \\
		\hfD &= \sum_{i=0}^p \hat{f}(\xi_i)l_i(\xi)\label{eq:fpoly}
	\end{align} 
	Here $p$ is the polynomial order and $l_i(\xi)$ are Lagrange polynomial bases in the reference domain, defined as:
	\begin{equation}
		l_i(\xi) = \prod_{\substack{j=0\\ j\ne i}}^p \bigg(\frac{\xi-\xi_j}{\xi_i-\xi_j}\bigg)
	\end{equation}
	The use of a superscript $\delta$ here highlights that a variable is a localised polynomial fit. Furthermore, for Eq.~(\ref{eq:fpoly}) the values of $\hat{f}(\xi_i)$ will have been calculated from $\hat{u}(\xi_i)$ and hence there is no enforcement on the polynomial fit in  Eq.~(\ref{eq:fpoly}) that it should be continuous between elements. For this reason we denote the polynomial fit with an added superscript \emph{D}, symbolising that the polynomial for the flux at this stage may not be $C_0$ continuous. 
	
\begin{figure}[tbhp]
		\centering
		\includegraphics[width=0.5\linewidth]{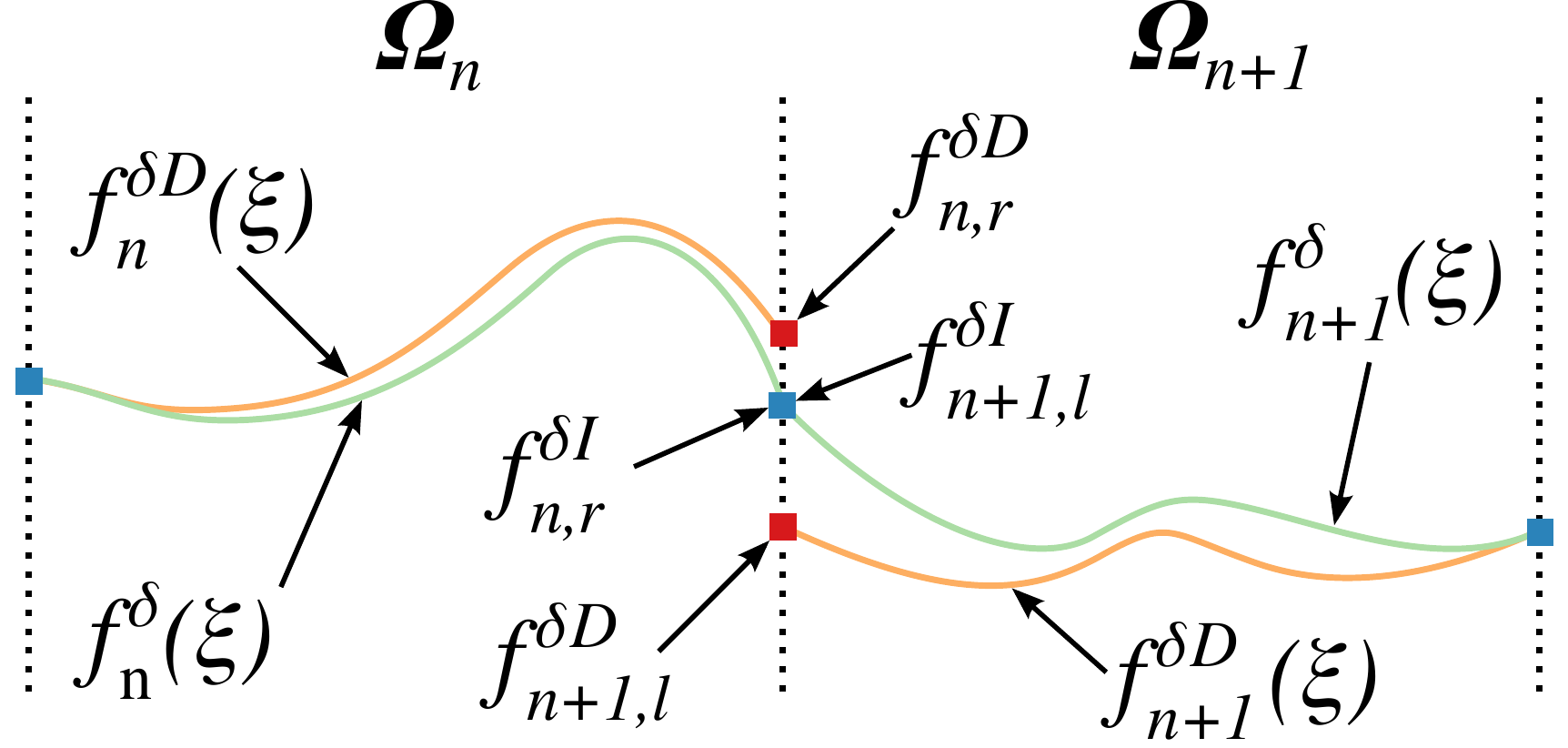}
		\caption{Diagram showing the notation and an interpretation of flux polynomial correction for one interface, cast into the physical domain for simplicity.}
		\label{fig:flux_interface}
	\end{figure}

	 The main component of the FR algorithm is the method by which $\hfD$ is made $C_0$ continuous. This is performed in a series of steps and a pictographic overview is given in Fig.~\ref{fig:flux_interface}. The first step is to interpolate the polynomial $\hud$ to the boundaries of the element. In 1D this is  defined as $\hat{u}^{\delta}(-1)=\hat{u}^{\delta}_l$  and $\hat{u}^{\delta}(1)=\hat{u}^{\delta}_r$. This process may be repeated for the discontinuous flux polynomial to obtain $\hfD_l$ and $\hfD_r$. 
	 
	 The second step is then to use collocated left and right interface values and fluxes to form a single common interface value for each element interface. There are several appropriate methods for achieving this, simple central differences can be used, however --- as is similarly experienced by other methods --- central differencing introduces no dissipation and so can be unstable without some form of stabilisation. For hyperbolic like equations, this problem can be solved by using a method that accounts for the upwind direction, therefore  introducing some dissipation. This includes methods such as that of Rusanov~\cite{Rusanov1961}, approximate Riemann solvers --- with Toro~\cite{Toro2009} providing a review of many methods --- or splitting methods. For the $j^\mathrm{th}$ element, we will denote the common interface fluxes as $\hfI_{j,l}$ and $\hfI_{j,r}$. Continuing, we will omit the $j$ and take it to mean this unless otherwise stated.
	 
 	The third step is to form the correction term that will force $\hfD$ to take the values of $\hfI$ at the interfaces and therefore enforce inter-element continuity. This procedure is performed by a left and right correction function, $h_l$ and $h_r$, such that the correction to the discontinuous flux is:
 	\begin{equation}
 		\hfC(\xi) = (\hfI_{l} - \hfD_{l})h_l(\xi) + (\hfI_{r} - \hfD_{r})h_r(\xi)
 	\end{equation}
 	The left and right correction functions hence have the properties that:
	\begin{align}
			h_l(-1) &=1, \quad h_l(1) = 0 \label{eq:hl_bc} \\
			h_r(-1) &=0, \quad h_r(1) = 1 \label{eq:hr_bc}
	\end{align}
	Therefore, the $C_0$ continuous flux polynomial can be written as:
	\begin{align}
		\hfd(\xi) &= \hfD(\xi) + \hfC(\xi) \\
		&= \hfD_{j}(\xi) + (\hfI_{l}-\hfD_{l})h_l(\xi) + (\hfI_{r}-\hfD_{r})h_r(\xi) \label{eq:correction_def}
	\end{align}
	
	Referring back to the equation to be solved, Eq.~(\ref{eq:linearConLaw}), we will now differentiate Eq.~(\ref{eq:correction_def}) for the $j^\mathrm{th}$ element within the numerical framework. This gives:
	\begin{align}
		\px{\hfd_j}{\xi} =& \px{\hfD_j}{\xi} + (\hfI_{j,l} - \hfD_{j,l})\dx{h_l}{\xi} + (\hfI_{j,r} - \hfD_{j,r})\dx{h_r}{\xi} \\
		=& \sum^p_{i=0}\hfD_{j,i}\dx{l_i(\xi)}{\xi} + (\hfI_{j,l} - \hfD_{j,l})\dx{h_l}{\xi} + (\hfI_{j,r} - \hfD_{j,r})\dx{h_r}{\xi}
	\end{align}
	
	This step results in the divergence of the correction function entering the solution, rather than the correction function itself. Therefore it is the divergence of the correction function that has to be fully contained in the polynomial space of FR and hence  $h_l$ and $h_r$ may be polynomials of order $p+1$. This can have the additional effect, through the extra information it provides, of raising the order of the corrected flux gradient to order $p$. Whereas, $\partial\hfD/\partial\xi$ is otherwise a polynomial of order $p-1$.
	
	Finally we may write:
	\begin{align}
		\px{\hud_j}{t} &= -\px{\hfd_j}{\xi} \\
	    & = -\sum_{i=0}^p\hfd_{j,i}\dx{l_i(\xi)}{\xi} - (\hfI_{j,l}-\hfd_{j,l})\dx{h_l(\xi)}{\xi} - (\hfI_{j,r}-\hfd_{j,r})\dx{h_r(\xi)}{\xi} \label{eq:FR_broken}
	\end{align}	
	At this stage a temporal integration method can be used to advance the solution in time. The result of this integration will be that recalculated flux polynomial may not be continuous again, and hence the correction procedure has to be repeated. 
	
	This presentation of the FR technique for first order 1D conservation equations should make it clear that the definition of the correction function is important to the exact nature of the method. In an early definition of the FR method, Huynh~\cite{Huynh2007} presented several correction functions and showed that their numerical characteristics could be quite different. For the remainder of this section we will detail works that have contributed to the definition of sets of correction functions. 
	
\begin{figure}[tbhp]
		\centering
		\includegraphics[width=0.4\linewidth]{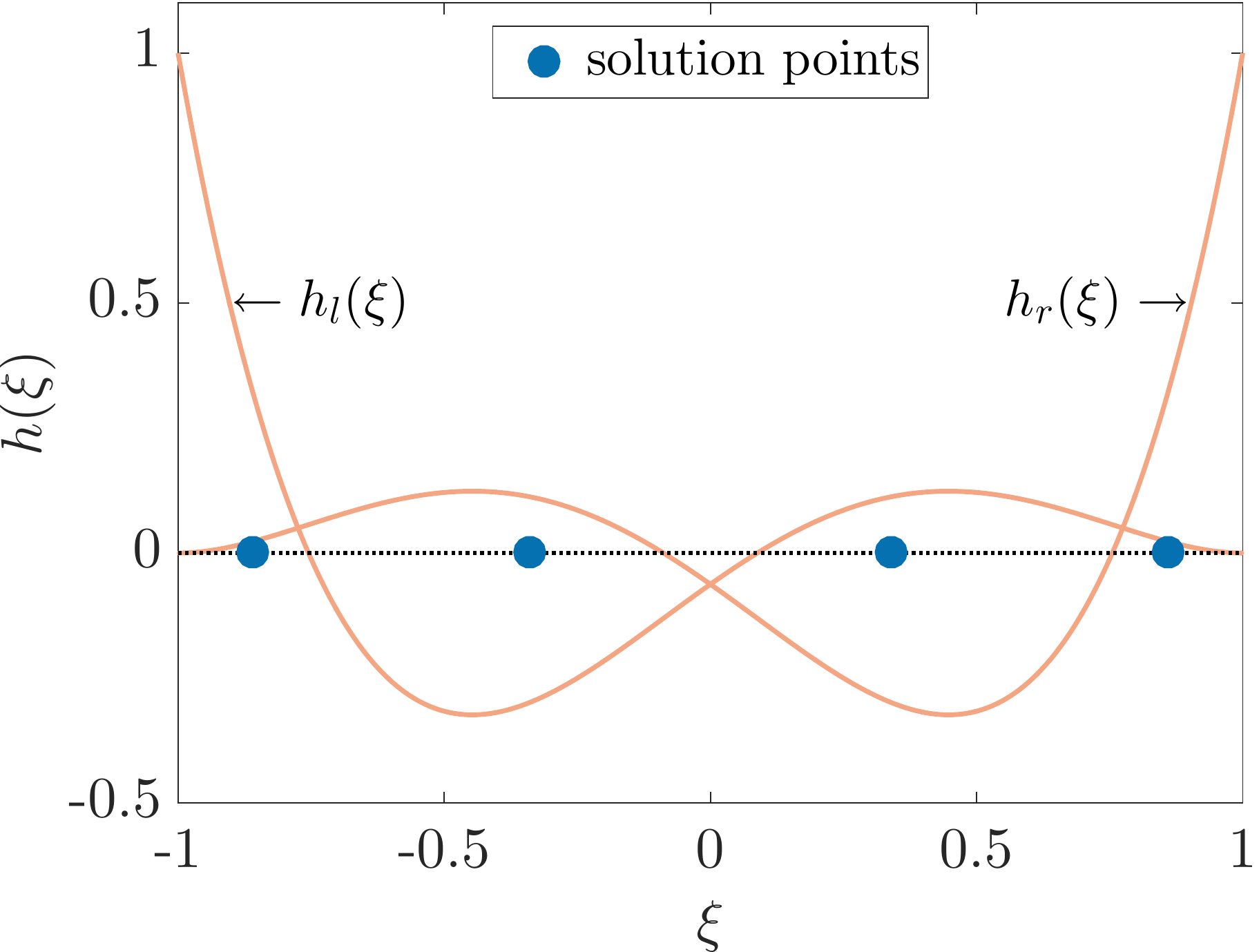}
		\caption{Left and right correction functions for $p=3$, $\iota=4/4725$, equivalent to Huynh's $g_2$ correction function~\cite{Huynh2007,Vincent2010}.}
		\label{fig:hu_corr}
	\end{figure}

	The first unifying set of stable correction functions was  proposed by Vincent~\etal~\cite{Vincent2010}, and we shall henceforth call this Original Stable FR (OSFR). This set is defined by a single variable as well as the polynomial order. The correction functions were found to be:  
	\begin{align}
		h_l &= \frac{(-1)^p}{2}\bigg[ \psi_p - \Big(\frac{\eta_p\psi_{p-1} + \psi_{p+1}}{1+ \eta_p}\Big)\bigg] \label{eq:osfr_hl}\\
		h_r &= \frac{1}{2}\bigg[ \psi_p + \Big(\frac{\eta_p\psi_{p-1} + \psi_{p+1}}{1+ \eta_p}\Big)\bigg] \label{eq:osfr_hr}\\
	\end{align}
	where:
	\begin{align}
		\eta_p &= \iota(2p+1)(a_pp!)^2 \label{eq:osfr_eta}\\
		a_p &= \frac{(2p)!}{2^p(p!)^2} 
	\end{align}
	and $\psi_i$ is the $i^{\mathrm{th}}$ order Legendre polynomials of the first kind defined on $\xi\in[-1,1]$, with $\iota$ taken as a free variable that we have \emph{modified to remove a factor of a half}. This set included some previously defined correction functions, such as Huynh's~\cite{Huynh2007} $g_2$ scheme (Fig.~\ref{fig:hu_corr}), the variation of SD that Jameson~\cite{Jameson2010} defined in FR, and the nodal Discontinuous Galerkin (DG) method.
	 
	Later an extension was made by Vincent~\etal~\cite{Vincent2015} to give what we will henceforth call Extended Range Stable FR (ESFR). This set of correction functions was far wider than that of OSFR, with more control variables. To introduce this set of correction functions, we will define the gradient of the left correction functions to be:
	\begin{equation}\label{eq:dcorr_leg}
		\dx{h_l}{\xi} = g_l(\xi) = \sum^p_{i=0} \tilde{\mathbf{g_l}}_{i}\psi_i(\xi)
	\end{equation}
	The right correction function is similarly defined. Here, and further in this text, we will use $\tilde{\bullet}$ to mean  something in a Legendre polynomial basis. It followed that the extended range of correction functions was defined by the equations:
	\begin{align}
		\tilde{\mathbf{g_l}} &= -\big(\tilde{\mathbf{M}} + \tilde{\mathbf{K}}\big)^{-1}\tilde{\mathbf{l}} \\
		\tilde{\mathbf{g_r}} &= \big(\tilde{\mathbf{M}} + \tilde{\mathbf{K}}\big)^{-1}\tilde{\mathbf{r}}
	\end{align}
	where $\tilde{\mathbf{l}} = [\psi_0(-1),\dots,\psi_p(-1)]^T$, $\tilde{\mathbf{r}} = [\psi_0(1),\dots,\psi_p(1)]^T$ and $\tb{M}$ is the Legendre polynomial mass matrix, defined as:
	\begin{align}
		\tilde{\mathbf{M}}_{i,j} = \int_{-1}^1 \psi_i\psi_j d\xi = \frac{2}{2j+1}\delta_{i,j}, \quad i,j \in \{0,\dots,p\}
	\end{align}
	
	The matrix $\tb{K}$ is freely set and dictates the form of the correction function. This matrix is constrained to be a real symmetric matrix and bound by:
	\begin{align}
		\tb{K} &= \tb{K}^{T} \\
		\tb{K}\tb{D} + (\tb{K}\tb{D})^T &= \mathbf{0} \\
		\tb{M} + \tb{K} &> \mathbf{0}
	\end{align}
	where
	\begin{align}
		\mathbf{D}_{i,j} &= \dx{l_j(\xi_i)}{\xi} \label{eq:Dmat} \\
		\mathbf{V}_{i,j} &= \psi_j(\xi_i) \\
		\tilde{\mathbf{D}} &= \mathbf{V}^{-1}\mathbf{D}\mathbf{V}
	\end{align}
	These conditions, together with $h_l(-1) = h_r(1) = 1$, $h_l(1) = h_r(-1) = 0$, and symmetry, defined a multi-parameter set of correction functions. Further, it was demonstrated that the single parameter OSFR set of Eq.(\ref{eq:osfr_hl}~\&~\ref{eq:osfr_hr}) is a subset of ESFR. To illustrate the correction function set, consider the case of $p=3$ from Vincent~\etal~\cite{Vincent2015}. The Legendre coefficients of the correction function gradient were found to be:
	\begin{equation}\label{eq:ESFR3}
		\tilde{\mathbf{g_l}} = -\begin{bmatrix}
			\frac{1}{2} \\[6pt]
			\frac{3(21\kappa_0 + 35\kappa_1 + 6)}{\upsilon} \\[6pt]
			\frac{5}{5\kappa_1 + 2} \\[6pt]
			\frac{21(5\kappa_1 + 2)}{\upsilon}
		\end{bmatrix} \quad \mathrm{where} \quad \upsilon = 175\kappa_1^2-42\kappa_0 - 12
	\end{equation}
	Here $\kappa_i$ are a series of variables that define the shape of the correction function.
\section{Generalised Sobolev Stability}
\label{sec:gsobolev}

	In Section~\ref{sec:fr} we presented the Flux Reconstruction method, from which it can be understood that an approximate solution to a conservative equation is formed in two parts. The derivative of a polynomial formed by fitting through localised data. Then secondly, sharing information between elements to ensure there is some degree of continuity between the local flux polynomials. This second step is achieved via a correction function. 
	
	The approximate solution that this process gives will inhabit some linear functional space. In particular a Sobolev space, which is a normed vector space, the norm of which will be described shortly. An important result in the development of FR was the recognition that the Sobolev space, and its norm, could be used to find a correction function that led to energy stability~\cite{Vincent2010,Jameson2010}. The result being that, in order for a scheme to be stable, it is required that the modified broken Sobolev norm exist and the energy contained within the broken Sobolev space decreases with time. Let us first define the broken Sobolev norm as:
	\begin{equation}\label{eq:sobolev_norm}
		\|u\|_{W^{p}_2} = \sum_{n=1}^m\|u\|_{n,W^{p}_2} = \sqrt{\sum^m_{n=1}\int_{\Omega_n}\sum_{i=0}^{p}(u^{(i)})^2\mathrm{d}\xi}
	\end{equation}
	where $u^{(i)}$ is the $i^{\mathrm{th}}$ spatial derivative of $u$ and $W^{p}_2$ is the $p^{\mathrm{th}}$ order $L_2$ Sobolev space. (In this case $W^{p}_2 = H^p$, where $H$ is a Hilbert space). If we consider linear advection with unit advective velocity ($f = u$), then a sufficient condition for energy stability is:
	\begin{equation}\label{eq:time_stab}
		\dx{}{t}\|u\|^2_{W^{p}_2} \leqslant 0.
	\end{equation}
	Although, from Hesthaven~\&~Warburton~\cite{Hesthaven2007}, this should be strictly defined to include the incoming and outgoing energy to form a necessary and sufficient condition. We will consider this later.
	
	In the derivation defining OSFR~\cite{Vincent2010}, the condition of Eq.~(\ref{eq:time_stab}) was satisfied, but through the use of a modified norm. The modified norm contained only the $i=0$ term and a weighted $i=p$ term, written as:
	\begin{equation}\label{eq:osfr_norm}
		\|u\|^2_{n,W^{p,\iota}_2} = \int_{\Omega_n}\sum_{i=0}^{p}(u)^2 + \iota(u^{(p)})^2\mathrm{d}\xi
	\end{equation}
	The reason for this is that the norm in Eq.~(\ref{eq:osfr_norm}), defines a metric and hence is sufficient to define the topology of the Sobolev space. With a Sobolev space defined, stability could be sought but with the complexity greatly reduced. However, this metric can be generalised to include all the derivative terms of the norm in Eq.~(\ref{eq:sobolev_norm}). The modified norm we propose is then:
	\begin{equation}\label{eq:sobo_mod}
		\|u\|^2_{n,W^{p,I}_2} = \int_{\Omega_n}\sum_{i=0}^{p}\iota_i(u^{(i)})^2d\xi
	\end{equation}
	where $\iota_i$ is a real constant and $I = \mathbf{I}_p = [\iota_0,\dots,\iota_p]^T$. By inspection, it can be seen that this follows the normal rules of norms and hence the metric space is complete in a Cauchy sense. This gives potential to find a greater set of metric spaces on the same topology, and hence a more general definition of correction functions --- which could have beneficial properties.
	
	To define a set of correction functions that satisfy this condition we will follow the method of Vincent~\etal~\cite{Vincent2010}, by using the flux reconstructed form of the first order conservation law to form the components of Eq.~(\ref{eq:sobo_mod}) applied to Eq.~(\ref{eq:time_stab}).  To perform this analysis we will then shift to consider the stability of a single element, therefore, from Hesthaven~\&~Warburton~\cite{Hesthaven2007}, the necessary and sufficient stability condition becomes:
	\begin{equation}\label{eq:sobo_mod_time}
		\dx{}{t}\|\hud\|^2_{W^{p,I}_2} = \dx{}{t}\rint{\sum^{p}_{i=0}\iota_i(\hudi{i})^2}{\xi} \leqslant -\sum^{p-1}_{i=0}\iota_i\bigg[\Big(\hudi{i}_r\Big)^2 - \Big(\hudi{i}_l\Big)^2 \bigg]
	\end{equation}
	Here the summation limit on the right-hand side has been taken to $p-1$, this is due to $\hud$ being a $p^\mathrm{th}$ order polynomial and, therefore, the higher-order terms cancel. If we set the flux function such that $\hfd = \hud$, then Eq.~(\ref{eq:FR_broken}) can be recast as:
	\begin{equation} \label{eq:fr_linear}
		\px{\hud}{t} = -\px{\hud}{\xi} - (\huI_l - \hud_l)\dx{h_l}{\xi} - (\huI_r - \hud_r)\dx{h_r}{\xi} 
	\end{equation}
	
	This equation can be used to form the individual terms required in Eq.~(\ref{eq:sobo_mod_time}). We will begin by forming the first derivative component by spatially differentiating Eq.~(\ref{eq:fr_linear}) and then multiplying by $\hudi{1}$ to give:
	\begin{multline}
		\px{\hud}{\xi}\dx{}{t}\bigg(\px{\hud}{\xi}\bigg) = \half\dx{}{t}\bigg(\px{\hud}{\xi}\bigg)^2 \\
		= -\px{\hud}{\xi}\pxi{2}{\hud}{\xi} - (\huI_l - \hud_l)\px{\hud}{\xi}\dxi{2}{h_l}{\xi} - (\huI_r - \hud_r)\px{\hud}{\xi}\dxi{2}{h_r}{\xi}
	\end{multline} 
	This can be simplified by further use of the product rule to:
	\begin{equation}\label{eq:fr_first}
		\half\dx{}{t}\bigg(\px{\hud}{\xi}\bigg)^2
		= -\half\px{}{\xi}\bigg(\px{\hud}{\xi}\bigg)^2 - (\huI_l - \hud_l)\px{\hud}{\xi}\dxi{2}{h_l}{\xi} - (\huI_r - \hud_r)\px{\hud}{\xi}\dxi{2}{h_r}{\xi}
	\end{equation} 
	
	Before continuing we should explain the nature of the differentiation that is being performed here, hence explaining the validity of applying the product rule. $\hud$ is being taking to mean a polynomial and, as it is continuous, so we mean differentiation to  be defined here using the standard limit approach. This links to a second point, that we are attempting to construct the components of the norm and not solve Eq.~(\ref{eq:fr_first}) in an FR framework and so each competent, say $\mathrm{d}(\hudi{1})^2/\mathrm{d}\xi$, is not required to lie in the same $p^\mathrm{th}$ order functional space. This is why differentiation is defined using a limit here and not using Eq.~(\ref{eq:Dmat}). This is the same method as was used by Vincent~\etal~\cite{Vincent2010}.
	
	Resuming our derivation, Eq.~(\ref{eq:fr_first}) can then be integrated over the reference domain and multiplied by $\iota_1~$ to give: 	
	\begin{multline}
		\half\dx{}{t}\rint{\iota_1(\hudi{1})^2}{\xi} \\
		= -\frac{\iota_1}{2}\rint{\px{}{\xi}\bigg(\px{\hud}{\xi}\bigg)^2}{\xi} - \iota_1(\huI_l - \hud_l)\rint{\px{\hud}{\xi}\dxi{2}{h_l}{\xi}}{\xi} \\
		- \iota_1(\huI_r - \hud_r)\rint{\px{\hud}{\xi}\dxi{2}{h_r}{\xi}}{\xi}
	\end{multline}
	This process of differentiating, multiplying by the derivative,	applying the product rule, and integrating can be repeated for all the derivatives up to the $p^\mathrm{th}$. The result is that we can write:
	\begin{multline}
		\half\dx{}{t}\rint{\iota_i(\hudi{i})^2}{\xi} = - \half\rint{\iota_i\px{}{\xi}\bigg(\pxi{i}{\hud}{\xi}\bigg)^2}{\xi} \\
		- (\huI_l-\hud_l)\rint{\iota_i\pxi{i}{\hud}{\xi}\dxi{i+1}{h_l}{\xi}}{\xi} \\ 
		- (\huI_r-\hud_r)\rint{\iota_i\pxi{i}{\hud}{\xi}\dxi{i+1}{h_r}{\xi}}{\xi}
	\end{multline}	
	This can then be substituted into Eq.~(\ref{eq:sobo_mod_time}), by summing over $i$ from $0$ to $p$, to give:
	\begin{multline}
		\half\dx{}{t}\|\hud\|^2_{W^{p,I}_2} = \half\dx{}{t}\rint{\sum^{p}_{i=0}\iota_i(\hudi{i})^2}{\xi} = - \half\underbrace{\rint{\sum^{p-1}_{i=0}\iota_i\px{}{\xi}\bigg(\pxi{i}{\hud}{\xi}\bigg)^2}{\xi}}_{I_d} \\
		- (\huI_l-\hud_l)\rint{\sum^p_{i=0}\iota_i\pxi{i}{\hud}{\xi}\dxi{i+1}{h_l}{\xi}}{\xi} \\
		- (\huI_r-\hud_r)\rint{\sum^p_{i=0}\iota_i\pxi{i}{\hud}{\xi}\dxi{i+1}{h_r}{\xi}}{\xi}	
	\end{multline}
	As $\hud$ is a $p^\mathrm{th}$ order polynomial, the term in $I_d$ for $i=p$ is identically zero and so the limits on the sum have been changed to reflect this.	Integration by parts then allows us to reformulate this as:
	\begin{multline}
		\half\dx{}{t}\|\hud\|^2_{W^{p,I}_2}  = - \half\rint{\sum^{p-1}_{i=0}\iota_i\px{}{\xi}\bigg(\pxi{i}{\hud}{\xi}\bigg)^2}{\xi} \\
		- (\huI_l-\hud_l)\sum^p_{i=0}\iota_i\Bigg[\bigg|\pxi{i}{\hud}{\xi}\dxi{i}{h_l}{\xi}\bigg|^1_{-1} -\rint{\pxi{i+1}{\hud}{\xi}\dxi{i}{h_l}{\xi}}{\xi} \Bigg] \\ 
		- (\huI_r-\hud_r)\sum^p_{i=0}\iota_i\Bigg[\bigg|\pxi{i}{\hud}{\xi}\dxi{i}{h_r}{\xi}\bigg|^1_{-1} -\rint{\pxi{i+1}{\hud}{\xi}\dxi{i}{h_r}{\xi}}{\xi} \Bigg]
	\end{multline} 
	If we apply the conditions that:
	\begin{align}
		\sum^p_{i=0}\iota_i\rint{\dxi{i}{h_l}{\xi}\pxi{i+1}{\hud}{\xi}}{\xi} &= \sum^p_{i=1}\iota_i\Bigg|\pxi{i}{\hud}{\xi}\dxi{i}{h_l}{\xi}\Bigg|^1_{-1} \label{eq:frgs_lcond}\\
		 \sum^p_{i=0}\iota_i\rint{\dxi{i}{h_r}{\xi}\pxi{i+1}{\hud}{\xi}}{\xi} &= \sum^p_{i=1}\iota_i\Bigg|\pxi{i}{\hud}{\xi}\dxi{i}{h_r}{\xi}\Bigg|^1_{-1} \label{eq:frgs_rcond}
	\end{align}
	then the remaining components are:
	\begin{multline}
		\half\dx{}{t}\|\hud\|^2_{W^{p,I}_2}  = - \half\rint{\sum^{p-1}_{i=0}\iota_i\px{}{\xi}\bigg(\pxi{i}{\hud}{\xi}\bigg)^2}{\xi} \\
		- \iota_0(\huI_l-\hud_l)\big|\hud h_l\big|^1_{-1} 
		- \iota_0(\huI_r-\hud_r)\big|\hud h_r\big|^1_{-1}
	\end{multline} 
	After the enforcement of the boundary conditions on $h_l$ and $h_r$ --- Eqs.~(\ref{eq:hl_bc}~\&~\ref{eq:hr_bc}) --- it can be written that:
	\begin{multline}
		\half\dx{}{t}\|\hud\|^2_{W^{p,I}_2}  = \half\iota_0\Big((\huI_l - \hud_l)\hud_l - (\huI_r - \hud_r)\hud_r  \Big)\\
		-  \half\sum^{p-1}_{i=0}\iota_i\Big[\big(\hudi{i}_r\big)^2 - \big(\hudi{i}_l\big)^2\Big]
	\end{multline}
	By then considering Eq.~(\ref{eq:sobo_mod_time}), we see that in order for the scheme to be stable the problem is reduced to proving:
	\begin{equation}
		\iota_0(\huI_l - \hud_l)\hud_l - \iota_0(\huI_r - \hud_r)\hud_r \leqslant 0
	\end{equation} 
 	This statement was explored by Vincent~\etal~\cite{Vincent2010} and found to be true for upwinded to centrally differenced interfaces.

\section{Spatial Stability Limits}
\label{sec:stability}

In order for the correction functions found using the above conditions to be valid, the norm defined by Eq.~(\ref{eq:sobo_mod}) must be positive and bound, \emph{i.e.} $0<\|\hud\|_{n,W^{p,I}_2}<\infty$. The other conditions for a valid norm will follow from differentiation being a linear operation. From the definition of the norm, thus: 
	\begin{equation}\label{eq:norm_limit}
		0 < \sqrt{\rint{\sum^p_{i=0}\iota_i(\hat{u}^{\delta,(i)})^2}{\xi}} < \infty  \quad \forall \; \hud \neq 0\\
	\end{equation}	
	The approximate solution, $\hud$, may then be cast into the Legendre polynomial basis as:
	\begin{equation}\label{eq:u_leg}
		\hud = \sum^p_{i=0}\tilde{u}_i\psi_i
	\end{equation}
	where $\psi_n$ is again a $n^\mathrm{th}$ order Legendre polynomial of the first kind. Upon substitution of Eq.~(\ref{eq:u_leg}) into Eq.~(\ref{eq:norm_limit}) we obtain:
	\begin{equation}
		0 < \sqrt{\rint{\sum^p_{i=0}\iota_i\Bigg(\sum^p_{j=0}\tilde{u}_j\dxi{i}{\psi_j}{\xi}\Bigg)^2}{\xi}} < \infty
	\end{equation}
	which may then be expanded to:
	\begin{equation}\label{eq:frgs_stability}
		0 < \sqrt{\sum^{p}_{j=0}\iota_0\bigg(\frac{2}{2j+1}\bigg)\tilde{u}_j^2 + \underbrace{\rint{\sum^{p-1}_{i=1}\iota_i\Bigg(\sum^p_{j=0}\tilde{u}_j\dxi{i}{\psi_j}{\xi}\Bigg)^2}{\xi}}_{I_m} + \iota_p\bigg(\frac{(2p)!}{2^pp!}\bigg)^2\tilde{u}_p^2 } < \infty
	\end{equation}
	
	A closed form for the integration of the product of two arbitrary Legendre polynomial derivatives, as in $I_m$, does exist and will be presented later. However, the exact evaluation of this condition will be dependent on the order, $p$, as well as on the solution, $\hud$. The effect the solution has on stability originates from the cross multiplication of $\tilde{u}_i$ terms, the origin of which can be understood by formulating the derivative of a Legendre polynomial as a Legendre series~\cite{Bailey1933,Dougall1952}. For example, there may be a $\tilde{u}_1\tilde{u}_3$ terms that arises in Eq.~(\ref{eq:frgs_stability}). Consequently, although Eq.~(\ref{eq:frgs_stability}) does constitute a necessary and sufficient condition for stability, a priori, we will only be able to form numeric limits on $\mathbf{I}_p$ that are sufficient conditions for stability due to the solution dependency.
	
	In the case of $\iota_i=0$ for $0<i<p$, i.e. when the scheme becomes OSFR, these cross multiplication terms vanish and the stability condition becomes that presented by Vincent~\etal~\cite{Vincent2010}. In this case the condition is formally necessary and sufficient for stability.
\section{Generalised Sobolev Correction Functions}
\label{sec:GSFR}

	Section~\ref{sec:gsobolev} presented analysis which showed that, subject to two conditions --- Eqs.~(\ref{eq:frgs_lcond}~\&~\ref{eq:frgs_rcond}) --- correction functions could be defined that generalised the Sobolev definition of energy stability in time. We are now concerned with finding the exact form of correction functions. The first step will be to express the left correction function in an orthogonal basis by using Legendre polynomials:
	\begin{equation}\label{eq:hl_leg}
		h_l(\xi) = \sum_{i=0}^{p+1} \tilde{\mathbf{h_l}}_{i}\psi_i(\xi) 
	\end{equation}
	Legendre polynomials are chosen as their weighting function in the orthogonality condition is unity, greatly simplifying later derivations. Therefore, substitution of Eq.(\ref{eq:u_leg}~\&~\ref{eq:hl_leg}) into Eq.(\ref{eq:frgs_lcond}) gives the condition on the left correction function being:
	\begin{equation}\label{eq:leg_condl}
	   \sum^p_{i=0}\iota_i\rint{\Bigg[\sum^{p+1}_{n=0}\sum^p_{m=0} \tilde{\mathbf{h_l}}_n\tilde{u}_m \dxi{i}{\psi_n}{\xi} \dxi{i+1}{\psi_m}{\xi}\Bigg]}{\xi} - \sum^p_{i=1}\iota_i\Bigg|\sum^{p+1}_{n=0}\sum^p_{m=0}\tilde{\mathbf{h_l}}_n\tilde{u}_m\dxi{i}{\psi_n}{\xi}\dxi{i}{\psi_m}{\xi}\Bigg|^1_{-1} = 0 
	\end{equation}
	Which may then be cast into a matrix form, with each row being for a value of $\tilde{u}_i$.
	\begin{equation} \label{eq:corr_pdef}
		\mathbf{L}_p\tilde{\mathbf{h_l}} = \begin{bmatrix} 0 & \hdots & 0 & 1 \end{bmatrix}^T = \mathbf{b_l}
	\end{equation} 	
	where the penultimate two entries of $\mathbf{L}_p$ and $\mathbf{b_l}$ are due to the enforcement of the boundary conditions on $h_l$. Hence, because of these boundary conditions, $\mathbf{L}_p$ is a square matrix. For the other rows, the entries of $\mathbf{L}_p$ are:
	\begin{equation}\label{eq:L_entry}
		\mathbf{L}_p[m-1][n] = \sum^p_{i=0}\iota_i\int^1_{-1}\dxi{i}{\psi_n}{\xi}\dxi{i+1}{\psi_m}{\xi} d\xi - \sum^p_{i=1}\iota_i\Bigg|\dxi{i}{\psi_n}{\xi}\dxi{i}{\psi_m}{\xi}\Bigg|^1_{-1}
	\end{equation}
	where $m,n\in\mathbb{N}$ with $1\leqslant m \leqslant p$ and $0\leqslant n \leqslant p+1$. The $m=0$ case is removed as it is identically zero, and, again, the final two rows will come from the boundary conditions on $h_l(\xi)$.
	
	For the right correction function, we may similarly write: 	
	\begin{equation}\label{eq:Lr}
		\mathbf{L}_p\tb{h_r} = \begin{bmatrix} 0 & \hdots & 1 & 0 \end{bmatrix}^T = \mathbf{b_r}
	\end{equation}
	Due to the important roll of correction functions in the definition of the scheme properties it is advantageous, but not essential, that the left and right correction functions are symmetric. We will now show that the GSFR correction function family is symmetric, allowing us to just focus on the left correction function. Symmetry implies the condition that:
	\begin{equation}\label{eq:lr_symm}
		\tb{h_l} = \mathrm{diag}(-1,1,-1\dots)\tb{h_r} = \mathbf{\Lambda}\tb{h_r}.
	\end{equation}	 	
	Hence, substituting Eq.~(\ref{eq:corr_pdef}), we may write:
	\begin{equation}
		\mathbf{L}_p\mathbf{\Lambda}\tb{h_r} = \mathbf{b_l}
	\end{equation}
	As all but the final two values of $\mathbf{b}$ are always zero it can be seen that solving this is equivalent to solving Eq.~(\ref{eq:Lr}). Therefore, GSFR defines a symmetric set of correction functions. Due to this, we will focus on the left correction function, with the right straightforwardly recovered through Eq.~(\ref{eq:lr_symm}).
	
	To evaluate Eq.~(\ref{eq:L_entry}), it can be useful to consider the results of Miller~\cite{Miller1963}. The results presented by Miller~\cite{Miller1963} were originally derived for associated Legendre polynomials of the first kind, but can be modified to consider just Legendre polynomials and the prescribed relationship of the derivatives, hence:
	\begin{equation}\label{eq:dmlegendre_int}
		\rint{\dxi{m}{\psi_n}{\xi}\dxi{m+1}{\psi_k}{\xi}}{\xi} = \sum^{\floor{\frac{n-m}{2}}}_{i=0}\sum^{\floor{\frac{k-m-1}{2}}}_{j=0}\frac{b_i(m,n)b_j(m+1,k)}{n+k-2(m+i+j)}[1-(-1)^{n+k-2(m+i+j)}]
	\end{equation}
	where we define:
	\begin{equation}\label{eq:legendre_b}
		b_i(m,n) = \frac{(-1)^i(2(n-i))!}{2^n(n-m-2i)!(n-i)!i!}
	\end{equation}
	Then for the gradient of Legendre polynomials at the end point:
	\begin{equation}
		\dxi{n}{\psi_j(-1)}{\xi} = \frac{(-1)^{j-n}(j+n)!}{2^nn!(j-n)!} \quad \mathrm{and} \quad \dxi{n}{\psi_j(1)}{\xi} = \frac{(1)^{j-n}(j+n)!}{2^nn!(j-n)!} \quad \mathrm{for} \quad j\geqslant n
	\end{equation} 
	which can be inferred from the work of Garfinkel~\cite{Garfinkel1964} and Holdeman~\cite{Holdeman1970}. Appendix~\ref{app:Lp} details pseudo-code for the generation of this matrix.  With these identities established, specific examples may now be evaluated.
\subsection{GSFR for \textit{p = 2}}
	For the case of $p=2$ the generalised correction function equation can be found to be given by:
	\begin{equation}\label{eq:corr_p2}
		\mathbf{L}_2\tilde{\mathbf{h_l}} = \begin{bmatrix}
			-\iota_0 & 0 & 3\iota_1 & 0 \\
			0 & -\iota_0 & 0 & 15(\iota_1+3\iota_2) \\
			1 & 1 & 1 & 1 \\
			1 & -1 & 1 & -1 
		\end{bmatrix} \tb{h_l} = \begin{bmatrix}
			0 \\ 0 \\ 0 \\ 1
		\end{bmatrix}
	\end{equation}
	And upon assessment of the limits presented by Eq.~(\ref{eq:frgs_stability}), the limits on $\mathbf{I}_2$ can be found to be:
	\begin{equation}
		\begin{bmatrix} 0 \\
			-\frac{1}{2}\Big(\frac{2}{3}\iota_0\Big) \\
			-\frac{1}{3}\Big(\frac{2}{5}\iota_0 + 6\iota_1\Big)
		\end{bmatrix} < \begin{bmatrix}
		\iota_0 \\
		\iota_1 \\
		\iota_2	\end{bmatrix}= 	
		\mathbf{I}_2 < \infty
	\end{equation}
\subsection{GSFR for \textit{p = 3}}
	\begin{equation}\label{eq:corr_p3}
		\begin{bmatrix}
			-\iota_0 & 0 & 3\iota_1 & 0 & 10\iota_1 \\
			0 & -\iota_0 & 0 & 15(\iota_1+3\iota_2) & 0 \\
			-\iota_0 & 0 & -(\iota_0-3\iota_1) & 0 & 15(3\iota_1+35\iota_2+105\iota_3) \\
			1 & 1 & 1 & 1 & 1 \\
			1 & -1 & 1 & -1 & 1
		\end{bmatrix} \tb{h_l} = \begin{bmatrix}
			0 \\ 0 \\ 0 \\ 0 \\ 1
		\end{bmatrix}
	\end{equation}
	Evaluating Eq.~(\ref{eq:frgs_stability}), the condition necessary for valid correction functions is: 
	\begin{multline}\label{eq:FRGS3_scond}
		0 < 2\iota_0\tilde{u}_0^2 + \bigg(\frac{2}{3}\iota_0 + \iota_1\bigg)\tilde{u}_1^2 + \bigg(\frac{2}{5}\iota_0 + 6\iota_1 + 18\iota_2\bigg)\tilde{u}_2^2 + \\
		\bigg(\frac{2}{7}\iota_0 + 8\iota_1 + 150\iota_2 + 255\iota_3\bigg)\tilde{u}_3^2 + \iota_1(\tilde{u}_1+2\tilde{u}_3)^2 < \infty
	\end{multline}
	Due to the final term in Eq.~(\ref{eq:FRGS3_scond}) being solution dependent, the transformation of this to a necessary condition on $\mathbf{I}_3$ is difficult. However, it can lead to a sufficient condition on the range of validity for $\mathbf{I}_3$:
	\begin{equation}
		\begin{bmatrix} 0 \\ 0\\
			-\frac{1}{18}\Big(\frac{2}{5}\iota_0 + 6\iota_1\Big) \\
			-\frac{1}{255}\Big(\frac{2}{7}\iota_0+8\iota_1+150\iota_2\Big)
		\end{bmatrix} < \mathbf{I}_3 < \infty	 
	\end{equation}

\subsection{GSFR for \textit{p = 4}}
	\begin{multline}\label{eq:corr_p4}
		\left[\begin{array}{cccccc}
			\iota_0 & 0 & 3\iota_1 &  0 \\
			0 & \iota_0 & 0 & 15(\iota_1+3\iota_2) \\
			\iota_0 & 0 & (\iota_0+3\iota_1) & 0 \\
			0 & \iota_0 & 0 & (-\iota_0 + 15\iota_1 +150\iota_2) & \\
			1 & 1 & 1 & 1 \\
			1 & -1 & 1 & -1 \\
		\end{array} \right. \\
		\left.\begin{array}{cc}
			10\iota_1 & 0 \\
			0 & (42\iota_1 + 315\iota_2) \\
			15(3\iota_1 + 35\iota_2 + 105\iota_3) & 0 \\
			0 & 105(\iota_1+31\iota_2 - 63\iota_3 + 945\iota_4) \\
		    1 & 1 \\
			1 & -1 \\
		\end{array}\right] \tb{h_l}  = \begin{bmatrix}
			0 \\ 0 \\ 0 \\ 0 \\ 0 \\ 1 \end{bmatrix}
	\end{multline}
	
	Evaluating Eq.~(\ref{eq:frgs_stability}) to find the necessary limits on validity: 
	\begin{multline}
		0 < 2\iota_0\tilde{u}_0^2 + \bigg(\frac{2}{3}\iota_0 + \iota_1\bigg)\tilde{u}_1^2 + \bigg(\frac{2}{5}\iota_0 + 2\iota_1 + 9\iota_2\bigg)\tilde{u}_2^2 + \bigg(\frac{2}{7}\iota_0 + 8\iota_1 + 150\iota_2 + 450\iota_3\bigg)\tilde{u}_3^2 \\
		 + \bigg(\frac{2}{9}\iota_0 + 11\iota_1 + 290\iota_2 + 7350\iota_3 + 11025\iota_4\bigg)\tilde{u}_4^2 \\ 
		 + \iota_1(2\tilde{u}_2+3\tilde{u}_4)^2 + \iota_1(\tilde{u}_1+2\tilde{u}_3)^2 + \iota_2(3\tilde{u}_2 + 20\tilde{u}_4)^2 < \infty
	\end{multline}
	And hence, with the same reasoning as for $p=3$, the sufficient conditions on $\mathbf{I}_4$ for valid correction functions are:
	\begin{equation}	
		\begin{bmatrix} 0 \\ 0 \\ 0 \\
			-\frac{1}{450}\Big(\frac{2}{7}\iota_0 + 8\iota_1 + 150\iota_2\Big) \\
			-\frac{1}{105^2}\Big(\frac{2}{9}\iota_0 + 11\iota_1 + 290\iota_2 + 7350\iota_3\Big)
		\end{bmatrix} < \mathbf{I}_4 < \infty
	\end{equation}
	
	At this point it is worth noting the recursive nature of the matrix $\mathbf{L}_p$. Hence the set of correction functions at $p$ is the union of $p^{\mathrm{th}}$ order correction functions and the correction functions defined by $\mathbf{L}_{p-1}$ and so on recursively down to the empty set. So if the $p^{\mathrm{th}}$ order correction function is defined by some eigenfunction $\mathbf{c}_p$ then the space of correction functions is defined as $\oplus^p_{i=0}\mathbf{c}_i$, \emph{i.e} each time the order is increased, one extra eigenfunction is introduced. However, in the special case of $\mathbf{I}_p = [1,0,\dots]^T$, the set of correction functions collapse to a single value, unique for each $p$. This is an interesting property that will be explored in future work.
	
	To show that these correction functions are in fact unique by comparison to both the OSFR and the ESFR, an attempt must be made to reconstruct the GSFR correction functions in both the ESFR and the  OSFR setting. Starting with OSFR, this method defines only the free parameter $\iota$, as described in Eq.~(\ref{eq:osfr_hl}~\&~\ref{eq:osfr_eta}). Taking $p=3$ and defining some arbitrary stable value of $\mathbf{I}_3$, the equivalent value of $\iota$ can be found using the value of $\tb{h_l}_{p+1}$ found from Eq.(\ref{eq:corr_p3}). Hence:
	\begin{equation}
		\iota = \frac{1}{(2p+1)(a_pp!)^2}\bigg(\frac{(-1)^{p+1}}{2\tb{h_l}_{p+1}}-1\bigg)
	\end{equation}
	
	To then compare to ESFR, it is easier to consider the gradient of $h_l$, defined in Eq.~(\ref{eq:dcorr_leg}), and defining the Legendre polynomial weights of GSFR similarly as $\tb{g_l}$. Then for the case of $p=3$ the corresponding ESFR weights can be found, using Eq.~(\ref{eq:ESFR3}) from Vincent~\etal~\cite{Vincent2015}, as:
	\begin{align}
		\kappa_1 &= - \bigg(\frac{1}{\tb{g_l}_2} + \frac{2}{5}\bigg)\label{eq:ESFR3_q1} \\
		\kappa_0 &= \frac{175\kappa_1^2\tb{g_l}_1+105\kappa_1-12\tb{g_l}_1+18}{42\tb{g_l}_1-63}
	\end{align}	
	It can then be further shown that to be an ESFR correction function for $p=3$, together with Eq.~(\ref{eq:ESFR3_q1}), the following must be satisfied:
	\begin{equation}
		\frac{175\kappa_1^2\tb{g_l}_3+105\kappa_1+42-12\tb{g_l}_3}{42\tb{g_l}_3} = \frac{175\kappa_1^2\tb{g_l}_1+105\kappa_1-12\tb{g_l}_1+18}{42\tb{g_l}_1-63}
	\end{equation}
	where $\kappa_1$ is defined by Eq.~(\ref{eq:ESFR3_q1}).
	
	\begin{figure}[tbhp]
	 	\centering
		\begin{subfigure}[b]{0.4\linewidth}
			\centering
			\includegraphics[width=\linewidth]{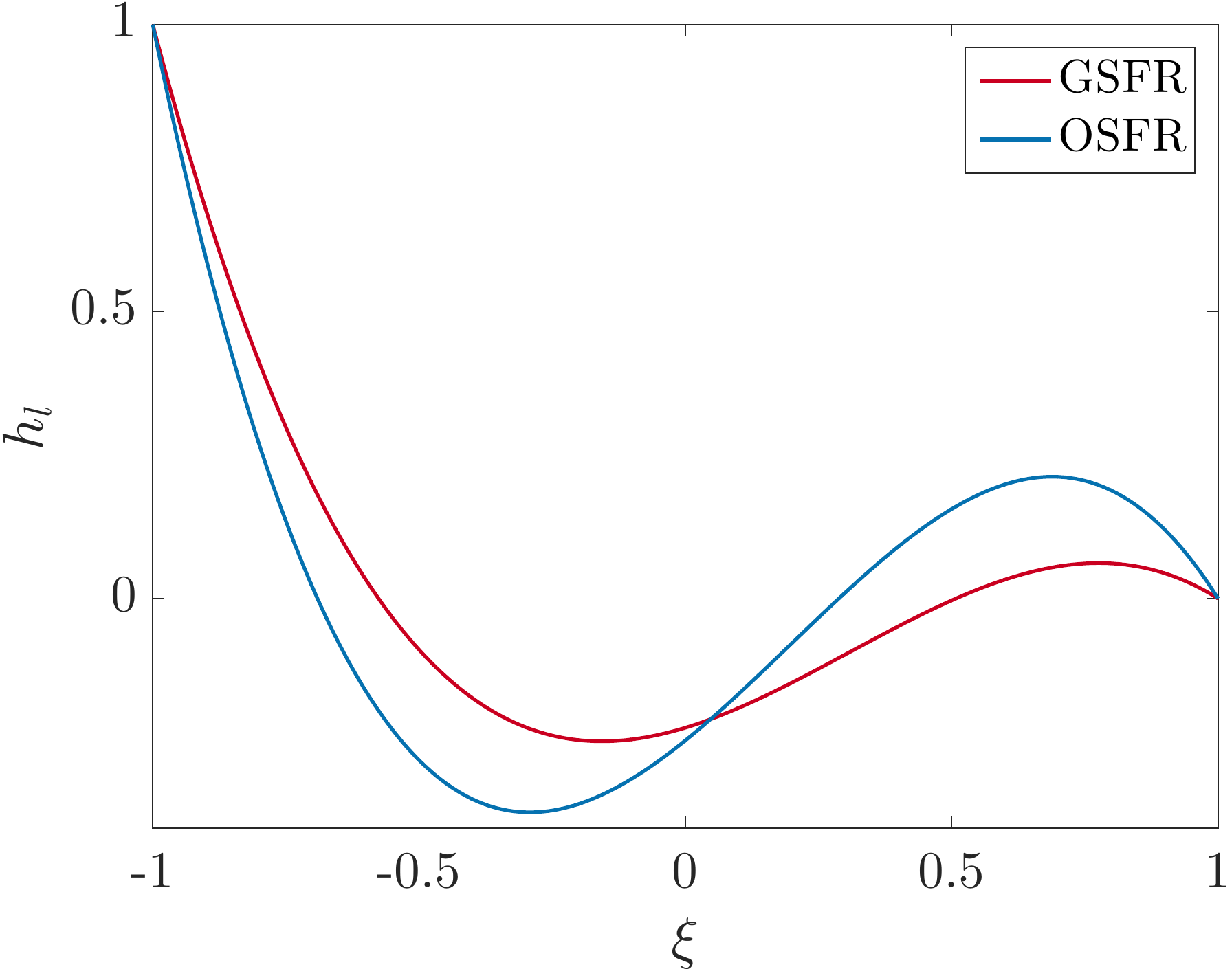}
			\caption{Left correction function comparison between OSFR and GSFR.}
			\label{fig:OSFR_GSFR_comp}
		\end{subfigure}
		~
		\begin{subfigure}[b]{0.4\linewidth}
			\centering
			\includegraphics[width=\linewidth]{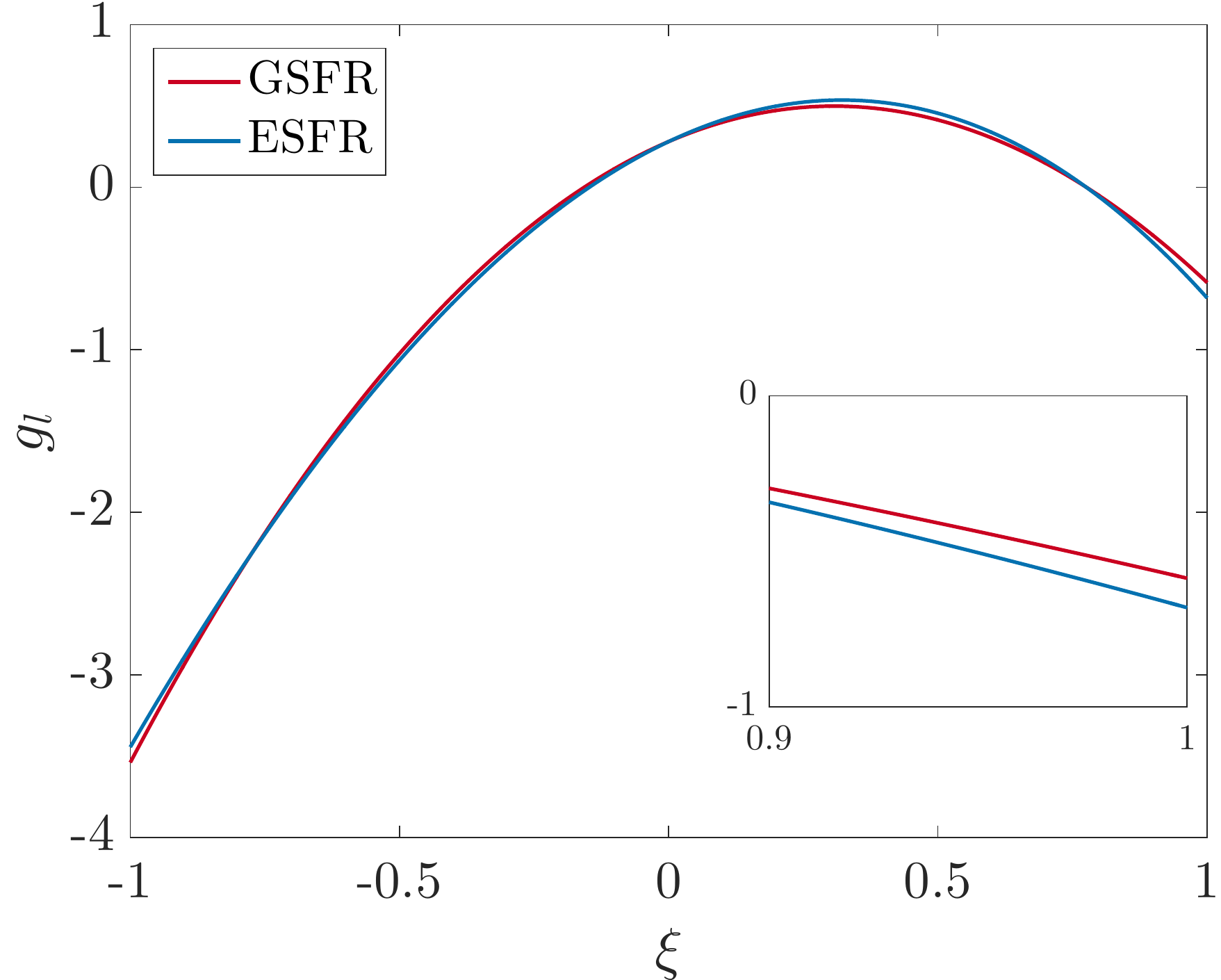}
			\caption{Left correction function gradient comparison between ESFR and GSFR.}
			\label{fig:ESFR_GSFR_comp}
		\end{subfigure}
		\caption{Comparison of OS, ES and GS correction functions. For $p=3$ and taking $\mathbf{I}_3= [1,\:0.01,\:0.01,\:0.1]^T$.}
		\label{fig:GSFR_uniquness}
	\end{figure} 		
	
	As is shown in Fig.~\ref{fig:GSFR_uniquness}, the correction functions found by satisfying Eq.~(\ref{eq:frgs_lcond}~\&~\ref{eq:frgs_rcond}) are in fact different from those defined by OSFR and ESFR. Furthermore, for OSFR it is trivial to show that it is a sub-set of GSFR, from there respective definitions, OSFR can be constructed when $\mathbf{I}_p = [1,\dots,\iota]^T$.  The ESFR  norm definition detailed by Vincent~\etal~\cite{Vincent2015} can be used to show ESFR is a subset of GSFR. This result could be expected as both OSFR and ESFR were found to be Sobolev stable. To find the corresponding values of $\mathbf{I}_p$ for a given $h_l$, originating from either OSFR or ESFR, the subject of the equation defining the GSFR correction must be change to $\mathbf{I}_p$. In the case of $p=3$ and setting $\iota_0=1$, this takes the form:
	\begin{equation}
		\underbrace{\begin{bmatrix}
			3\tb{h_l}_2 + 10\tb{h_l}_4 & 0 & 0 \\
			15\tb{h_l}_3 & 45\tb{h_l}_3 & 0 \\
			3\tb{h_l}_2+45\tb{h_l}_4 & 525\tb{h_l}_4 & 1575\tb{h_l}_4
		\end{bmatrix}}_{\mathbf{H}_3} \begin{bmatrix}
		\iota_1 \\ \iota_2 \\ \iota_3
		\end{bmatrix} = \begin{bmatrix}
			\tb{h_l}_0 \\ \tb{h_l}_1 \\ \tb{h_l}_0 + \tb{h_l}_2
		\end{bmatrix}
	\end{equation}
	Hence, ESFR and OSFR are both recoverable from GSFR. The exception to the invertibility of $\mathbf{H}_3$ over the set of ESFR and OSFR is when either $\tb{h_l}_3=0$ or $\tb{h_l}_4=0$, as GSFR has multiple ways of constructing these lower order correction functions.
\section{Fourier and von Neumann Analysis}
\label{sec:vn}

Initial characterisation of the set of schemes defined by the GSFR correction functions is investigated via von Neumann analysis applied to the linear advection equation. With unit wave speed, this can be cast as the semi-discrete matrix equation for flux reconstruction:
	\begin{equation}
		\frac{\partial \mathbf{u}_j}{\partial t} = - \Big(J_{j+1}^{-1}\mathbf{C}_{+1}\mathbf{u}_{j+1} + J_j^{-1}\mathbf{C}_0\mathbf{u}_j + J_{j-1}^{-1}\mathbf{C}_{-1}\mathbf{u}_{j-1}\Big)
	\end{equation}
	where $J_j$ is the $j^{\mathrm{th}}$ element's Jacobian and defining the operator matrices as:
	\begin{align}
		\mathbf{C}_{+1} &= (1-\alpha)\mathbf{g_r}\mathbf{l_l}^T \\
		\mathbf{C}_0 &= \mathbf{D} - \alpha\mathbf{g_l}\mathbf{l_l}^T - (1-\alpha)\mathbf{g_r}\mathbf{l_r}^T\\		
		\mathbf{C}_{-1} &= \alpha\mathbf{g_l}\mathbf{l_r}^T
	\end{align}
	where $\alpha$ is the upwinding ratio with $\alpha=1$ corresponds to upwinding and $\alpha=0.5$ to central difference. The array $\mathbf{g_l}$ is the value of $g_l(\xi)$ at the solution quadrature points and $\mathbf{l_l}$ is a vector interpolating from the solution points to the left interface. With $\mathbf{g_r}$ and $\mathbf{l_r}$ similarly defined. This method in this form was initially presented by Trojak~\etal~\cite{Trojak2019} and is similar to that of Huynh~\cite{Huynh2007}. For the study to be performed we will focus on the case when applied to a uniform mesh and hence $J$ will be constant. If a Bloch wave is then defined such that: 
	\begin{equation}
		u(x,t) = v \exp{\big(i(kx - \omega t)\big)}
	\end{equation}
	then the discretised equation can be found to be:
	\begin{equation}
		\frac{\partial \mathbf{u}_j}{\partial t} = -\Big(J_{j+1}^{-1}\mathbf{C}_{+1}\exp{\big(-ik\delta_{j}\big)} + J_j^{-1}\mathbf{C}_0 + J_{j-1}^{-1}\mathbf{C}_{-1}\exp{\big(-ik\delta_{j-1}\big)} \Big)\overline{\mathbf{u}}_j = \mathbf{Q}(k)\overline{\mathbf{u}}_j
	\end{equation}
	where $\delta_j = x_j - x_{j-1}$. Therefore, after differentiating the Bloch wave w.r.t. time, the wave phase velocity can be found from the eigenvalue problem:
	\begin{equation}
		c(k) \mathbf{v} = \frac{i}{k}\mathbf{Q}(k)\mathbf{v}
	\end{equation}
	
   Consequently, the dispersion and dissipation can be found as $\Re{(\hat{\omega})} = \Re{(c)}\hat{k}$ and $\Im{(\hat{\omega})} = \Im{(c)}\hat{k}$ respectively. To investigate the fully-discretised equation, a temporal discretisation can be included by defining the fully-discretised linear operator, or update matrix, as:
   	\begin{align}
		\mathbf{u}^{n+1}_j &= \mathbf{R}(\mathbf{Q})\mathbf{u}^n_j \\
		\mathbf{R}_{44} &= \sum^4_{n=0}\frac{(\tau\mathbf{Q})^n}{n!}
	\end{align}
	where $\tau$ is the explicit time integration step and $\mathbf{R}_{44}$ gives an example definition of the update matrix for RK44 temporal integration. Therefore, for von Neumann type spatial-temporal stability to be assured~\cite{Isaacson1994}, the spectral radius of the update matrix must be less than or equal to 1, $\rho(\mathbf{R}) \leqslant 1$.

	As was discussed in Section~\ref{sec:GSFR}, the cross product between some solution weights, for example the last term in Eq.~(\ref{eq:FRGS3_scond}), leads to difficulty in transforming the necessary condition for a positive definite norm to a condition on $\mathbf{I}_p$. Although a sufficient condition on $\mathbf{I}_p$ can be found, $\mathbf{I}_p$ will be allowed to vary beyond these bounds to aid understanding in the underlying necessary condition. Throughout this analysis $\iota_0$ will also be taken as one; it should be understood that this is an arbitrary choice, but that choosing a different value will simply lead to a linear scaling of all other $\iota_i$ values shown here. 
	
	\begin{figure}[tbhp]
	 	\centering
	 	\begin{subfigure}[b]{0.45\linewidth}
			\centering
			\includegraphics[width=\linewidth]{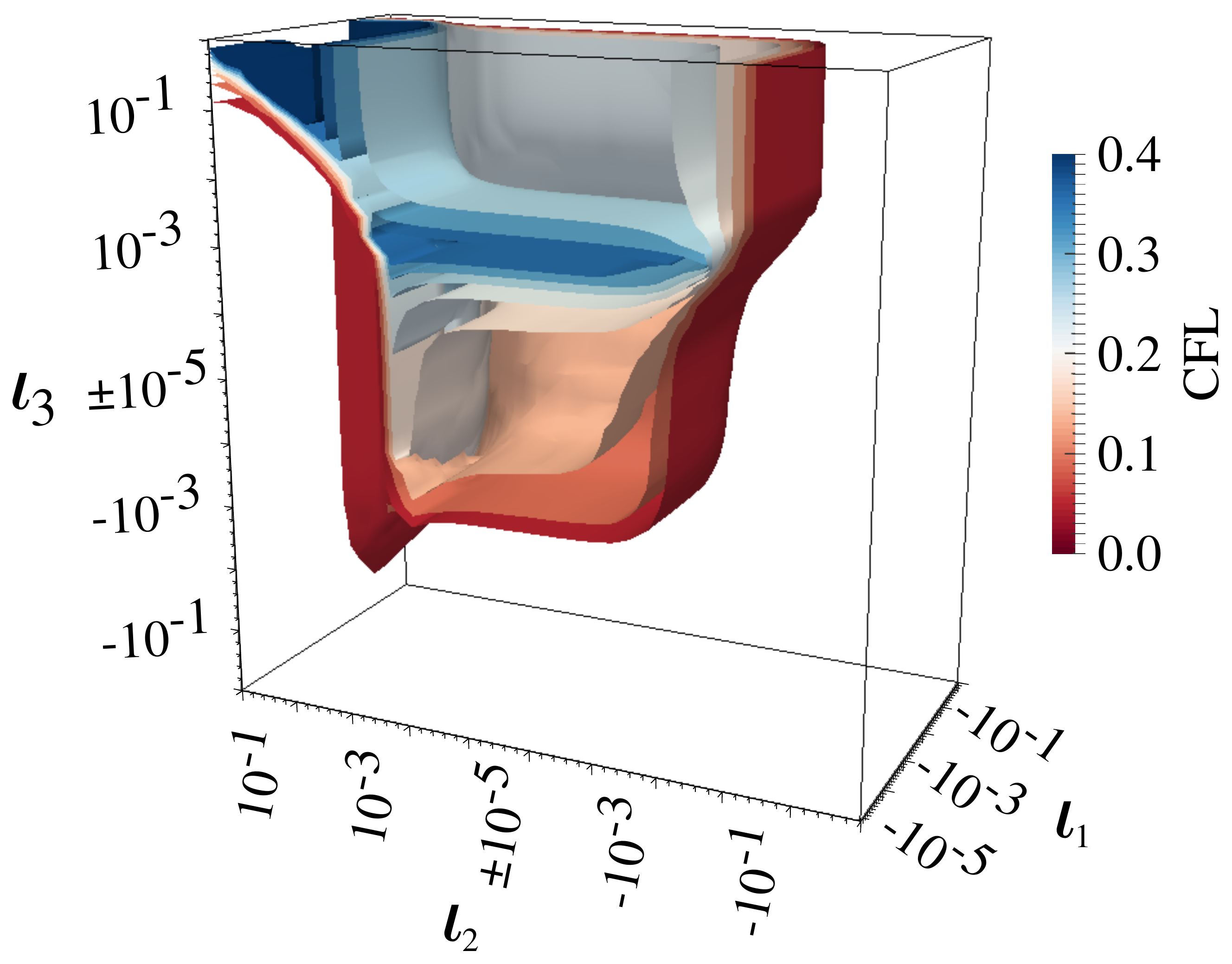}
			\caption{$\iota_1<0$}
			\label{fig:FRGS3_CFL_mi1}
		\end{subfigure}
		~
		\begin{subfigure}[b]{0.45\linewidth}
			\centering
			\includegraphics[width=\linewidth]{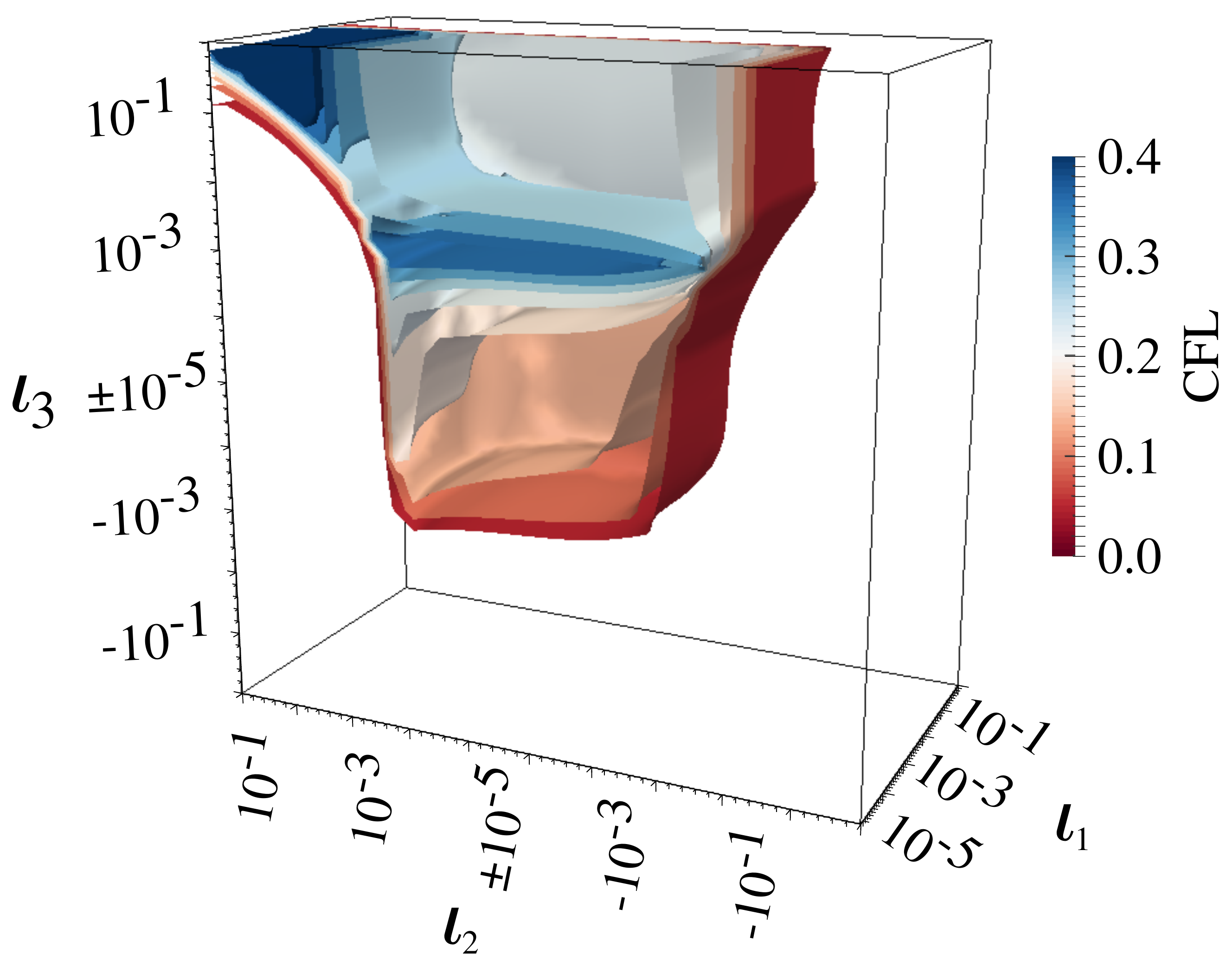}
			\caption{$\iota_1>0$}
			\label{fig:FRGS3_CFL_pi1}
		\end{subfigure}
		~
		\begin{subfigure}[b]{0.45\linewidth}
			\centering
			\includegraphics[width=\linewidth]{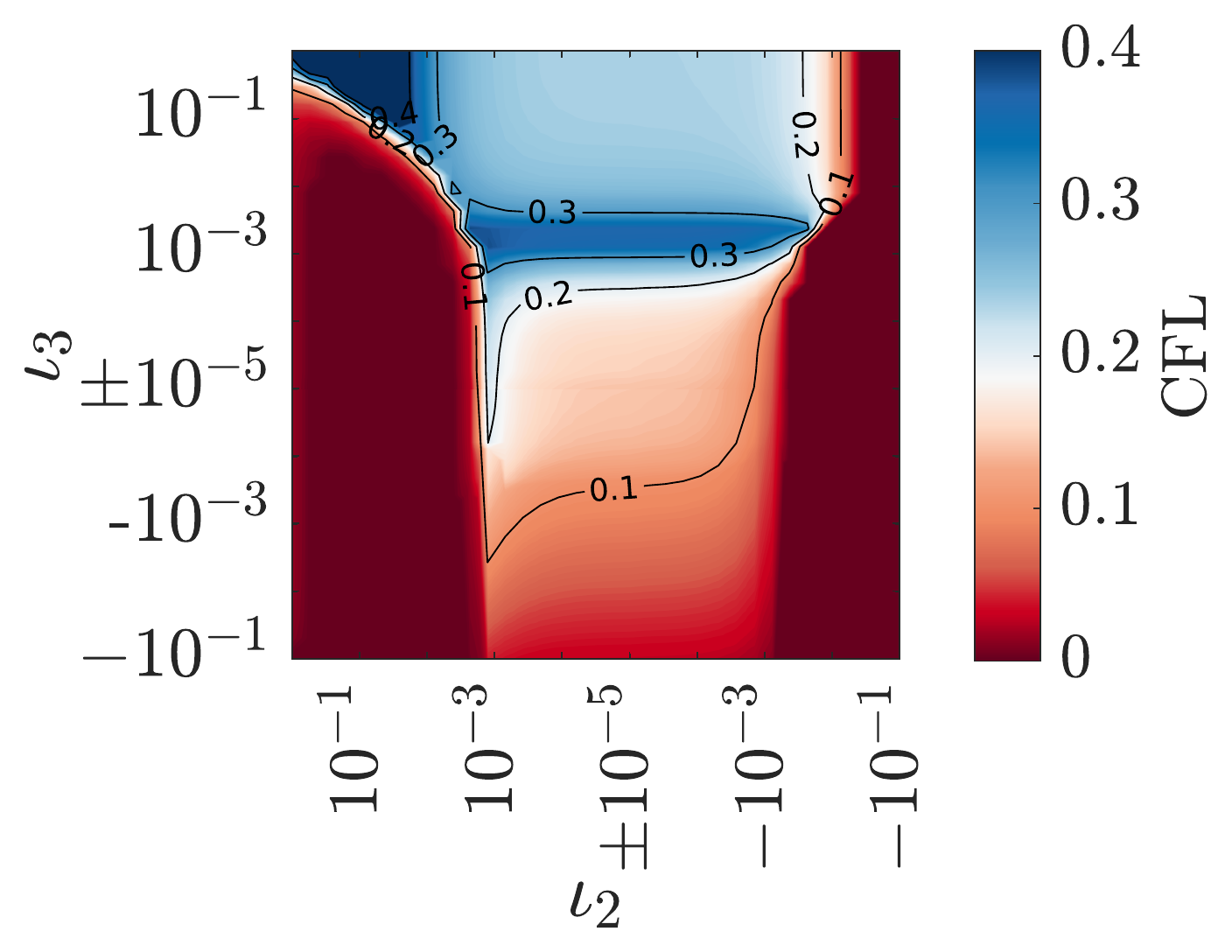}
			\caption{$\iota_1=0$}
			\label{fig:FRGS3_CFL_slice}
		\end{subfigure}
		\caption{CFL limit for upwinded FR with GSFR correction functions, $p=3$, and RK44 temporal integration on a regular grid.}
		\label{fig:FRGS3_CFL}
	\end{figure}
	
	Beginning with the case of $p=3$, Fig.~\ref{fig:FRGS3_CFL} shows the variation of CFL number with $\mathbf{I}_3$ for low storage RK44 temporal integration. In this figure positive and negative log axes in $\iota_2$ and $\iota_3$ are used to clearly show the CFL manifold, the join is at $\pm10^{-5}$. To clarify some key behaviour, Fig.~\ref{fig:FRGS3_CFL_slice} shows a slice for $\iota_1=0$. Through these figures, the extent of the stable domain of $\mathbf{I}_3$  can be seen. There are several interesting regions which can be seen, first of which is that as $\iota_3,\iota_2 \rightarrow \infty$. In this case the correction functions converges upon those constructed for $p=2$, which is similar to the behaviour observed from ESFR schemes. The second region of interest is where $\mathbf{I}_3 \approx [1,0,10^{-3},10^{-3}]^T$, here the CFL limit can be seen to have a local maximum. This region of local maximum is also believed to give recovery of high-order. This is due to the scale of the Legendre weights forming the correction and the position of the Nyquist wavenumber. It has previously been observed~\cite{Vincent2011} that a drop in order of accuracy reduces the Nyquist limit. To exactly find the order of accuracy~(OOA) a numerical method will be introduced in Section~\ref{sec:numerical}.
	
	\begin{figure}[tbhp]
	 	\centering
		\begin{subfigure}[b]{0.4\linewidth}
			\centering
			\includegraphics[width=\linewidth]{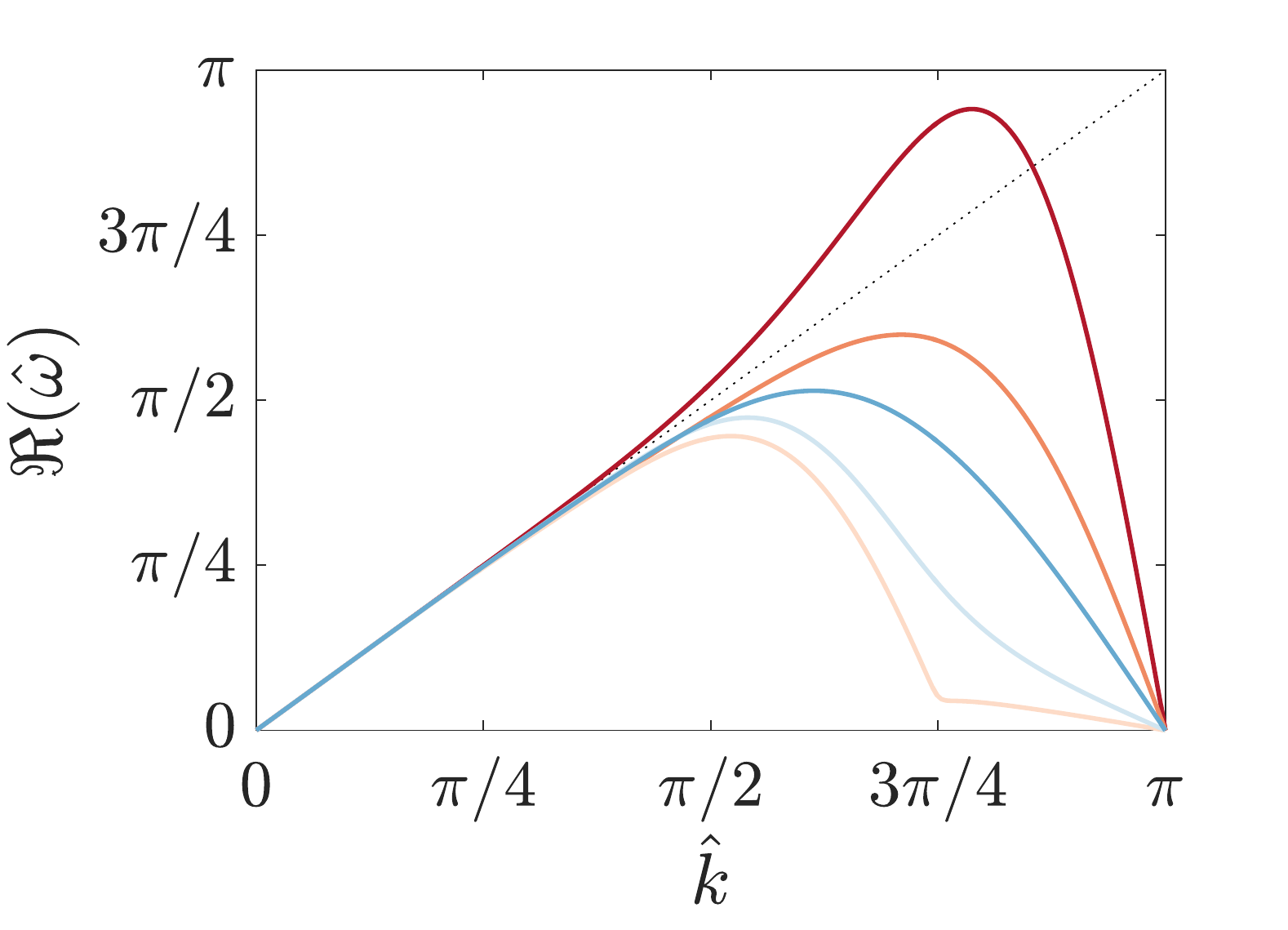}
			\caption{Dispersion.}
			\label{fig:FRGS3_comp_R}
		\end{subfigure}
		~
		\begin{subfigure}[b]{0.4\linewidth}
			\centering
			\includegraphics[width=\linewidth]{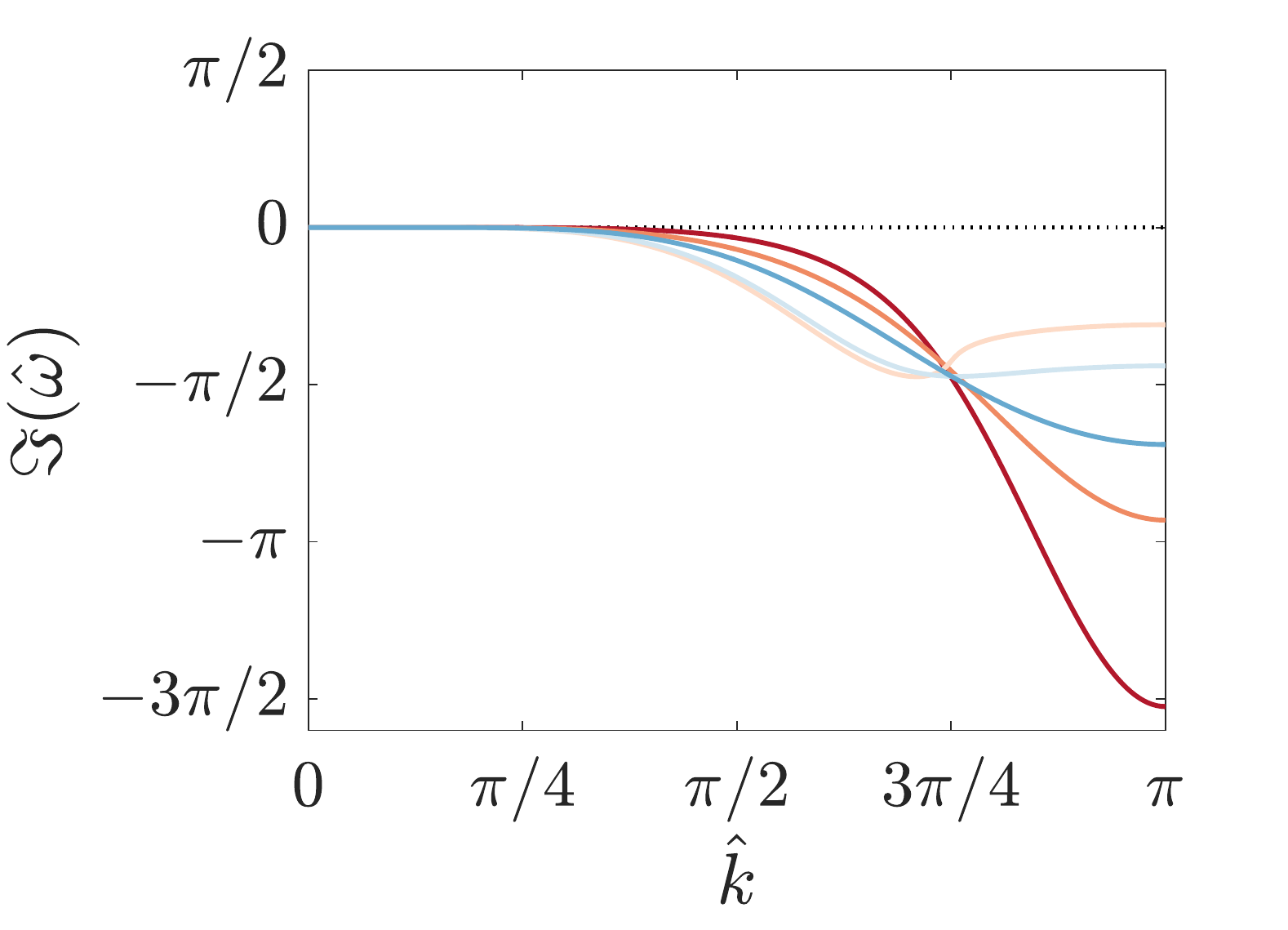}
			\caption{Dissipation.}
			\label{fig:FRGS3_comp_I}
		\end{subfigure}
		~
		\begin{subfigure}[b]{0.4\linewidth}
			\centering
			\includegraphics[width=\linewidth,trim= 0mm 30mm 0mm 25mm,clip=true]{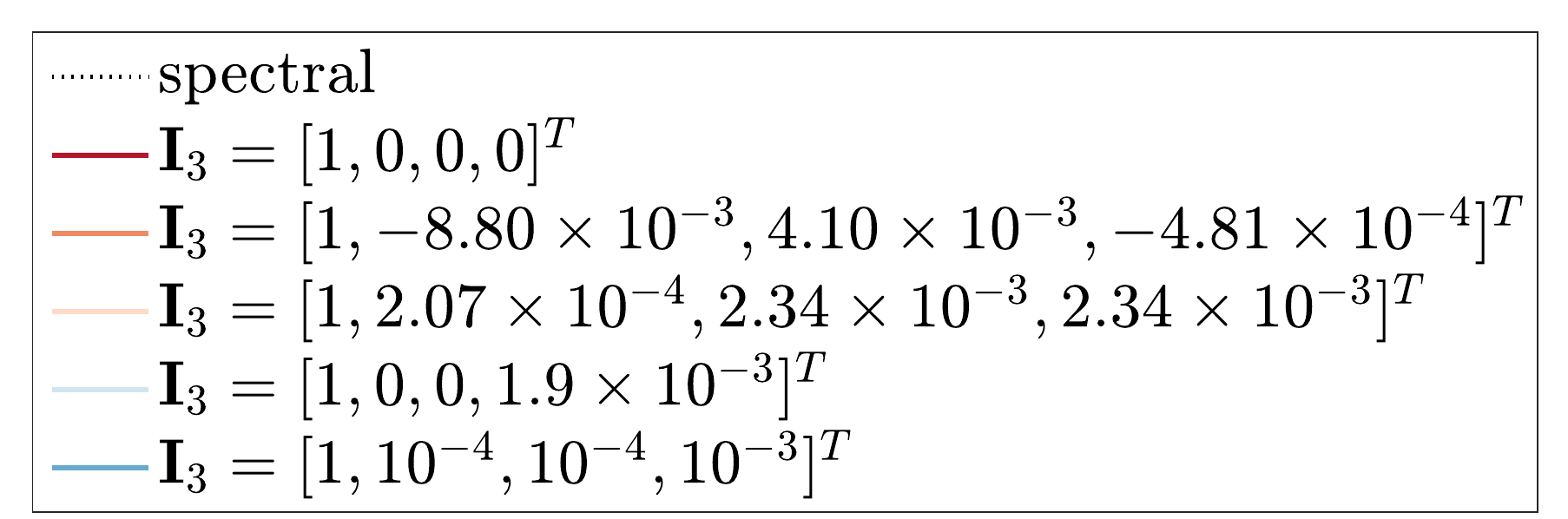}
		\end{subfigure}
		\caption{Dispersion and dissipation for selected $p=3$ GSFR correction functions with interface upwinding.}
		\label{fig:FRGS3_dispdiss}
	\end{figure} 
	
	The dispersion and dissipation characteristics for some selected correction functions are presented in Fig.~\ref{fig:FRGS3_dispdiss} for $p=3$.  The correction function tested are: DG; the ESFR correction function that gives optimal temporal stability~\cite{Vermeire2016}; the GSFR optimal correction function; the OSFR correction function found to give optimal temporal stability with this combination of spatial-temporal scheme~\cite{Vincent2011}; and lastly a stable GSFR correction function chosen arbitrarily. 
	
	For correction three, the Nyquist wavenumber is maintained at that for an OOA of four as it was previously suspected it may. However, there is a significant region where the phase velocity is low, i.e $c_p=\hat{\omega}/\hat{k} \approx 0$, which may cause additional dispersion errors. Additionally, at high wavenumbers, there is a significant reduction in the dissipation exhibited. This combination of dispersion and dissipation seems to have given rise to an increased CFL limit, but potentially at the cost of spectral performance and order accuracy, both of which will be investigated later.
	
	\begin{figure}[tbhp]
	 	\centering
		\begin{subfigure}[b]{0.45\linewidth}
			\centering
			\includegraphics[width=\linewidth]{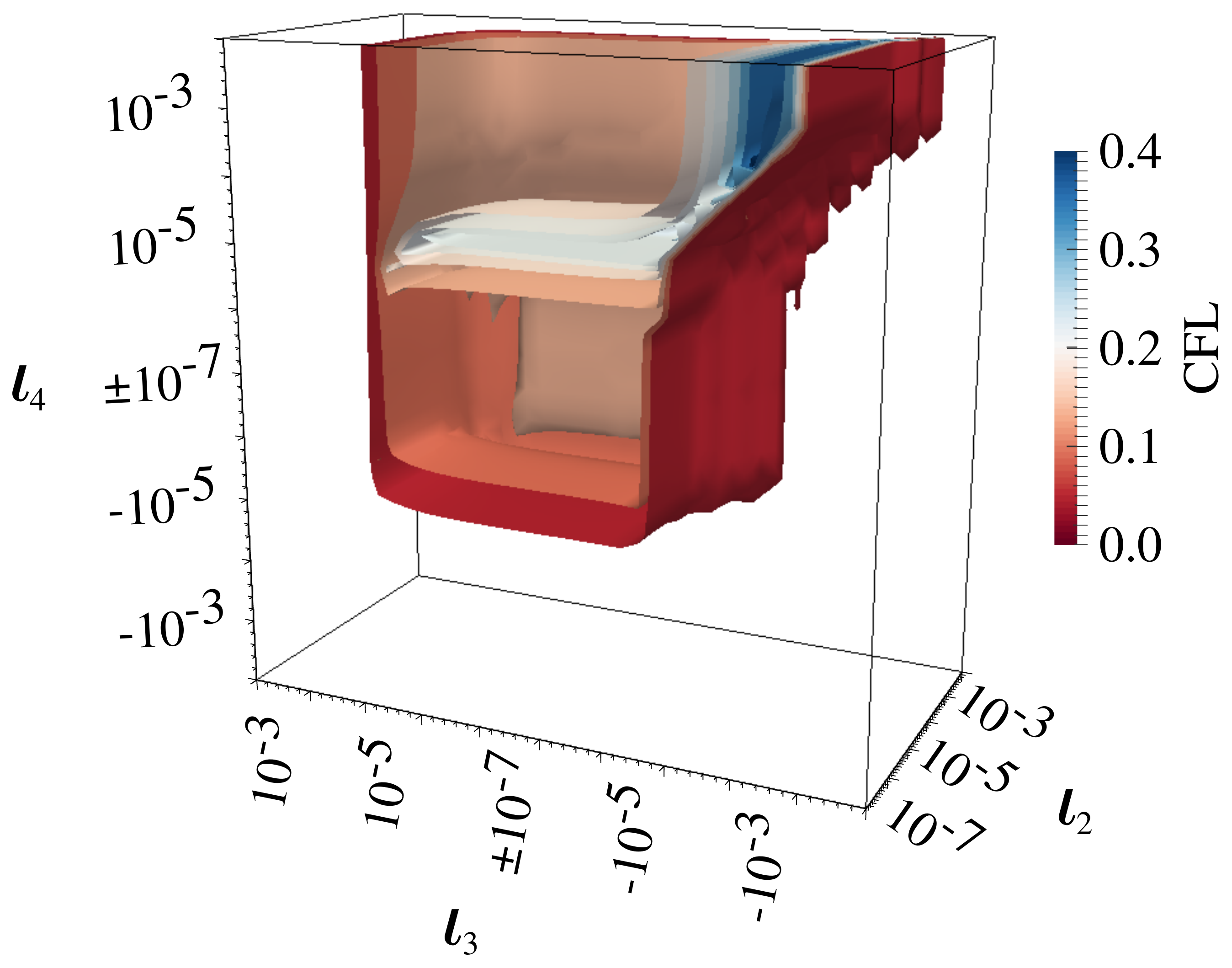}
			\caption{CFL limit isosurfaces.\\ \quad \quad}
			\label{fig:frgs4_cfl_iso}
		\end{subfigure}
		~
		\begin{subfigure}[b]{0.45\linewidth}
			\centering
			\includegraphics[width=\linewidth]{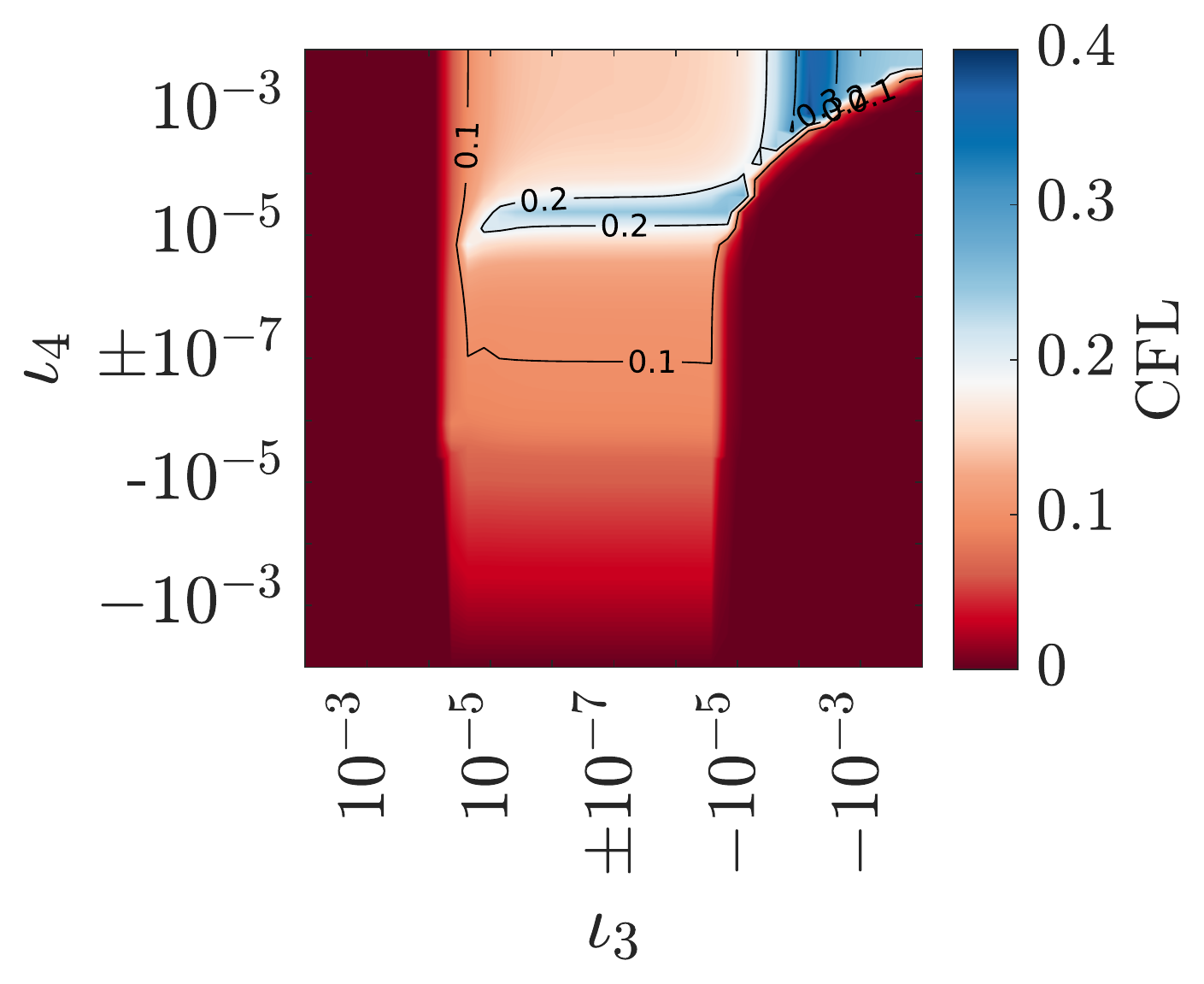}
			\caption{CFL limit contours for $\iota_1=10^{-5}$ and $\iota_2=10^{-7}$.}
			\label{fig:frgs4_cfl_con}
		\end{subfigure}
		\caption{CFL limit for upwinded FR with GSFR correction functions, $p=4$, and RK44 temporal integration on a regular grid. ($\iota_1 = 1\times 10^{-5}$).}
		\label{fig:FRGS4_CFL}
	\end{figure}	
	
	Returning to the investigation of CFL limits, Fig.~\ref{fig:FRGS4_CFL} shows two regions of high CFL limits which are similar to those shown for $p=3$. First considering the region shown in Fig.~\ref{fig:frgs4_cfl_con} at $\iota_3<-10^{-4}$ and $\iota_4>10^{-4}$. It is thought that this region contains the degradation to $p=3$. Figure~\ref{fig:frgs4_cfl_iso} clearly shows a localised maximum CFL limit, the free parameter is four dimensional, and hence a search method will be employed to find the global high-order recovering maximum CFL limit. This method, together with the results, will be introduced in Section~\ref{sec:numerical}.
\section{Numerical Tests}
\label{sec:numerical}

\subsection{Linear Homogeneous Advection}

Numerical tests were undertaken to validate analytical findings and provide useful further insight. The first such numerical test was contrived to calculate the OOA of the scheme while solving the linear advection equation with unit advection velocity. The domain set-up for this case was a periodic domain with $x\in[0,2\pi]$, with the number of elements allowed to vary such that several degrees of freedom could be tested. In particular, the numbers of elements used were $N = [50,55\dots 75]$. The interfaces were fully upwinded and a time step was chosen such that the temporal integration was not the primary source of error. The initial solution was taken to be a plane wave, with a wavenumber that would be expected to be well resolved on the grids:
	\begin{equation}
		u(x;t=0) = u^e(x;t=0) = \cos{(kx)}, \quad k=\frac{1}{2\pi}
	\end{equation}	
	where $u^e(x,t)$ is the exact solution and $u(x,t)$ is the approximate solution from FR. Taking the solution after some number of time integration steps, the point averaged error can be calculated for several grid densities, with $n_s$ points: 
	\begin{align}
		\epsilon_2 &= \frac{1}{n_s}\sum_{i=1}^{n_s} |u^e(x_i,t)-u(x_i,t)| \\
		\epsilon_2 &= |u^e(x,t) - u(x,t)| = \mathcal{O}(n_s^{-l}) \label{eq:ooa_def}
	\end{align}
	where $l$ is the numerically realised OOA, which can be found through taking logs of Eq.(\ref{eq:ooa_def}). With this established, the correction functions can now be varied to show the effect of $\mathbf{I}_p$ on OOA. 
	
	\begin{figure}[tbhp]
		\centering
		\includegraphics[width=0.6\linewidth]{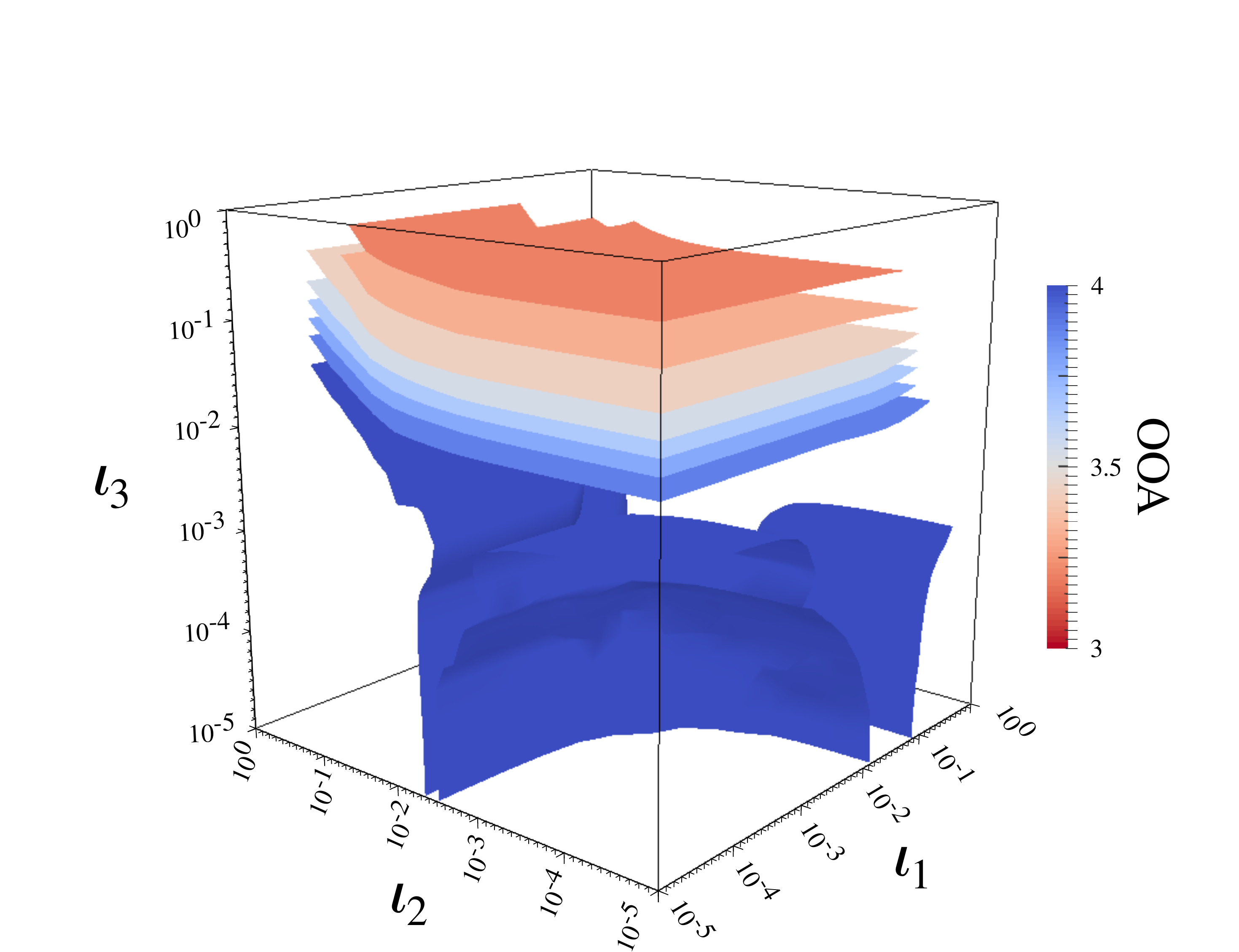}
		\caption{Order of accuracy of GSFR, $p=3$ for linear advection of a wave with $k=1/2\pi$, at $t=\pi$.}
		\label{fig:FRGS3_OOA}
	\end{figure}	
	
	Figure~\ref{fig:FRGS3_OOA} shows the OOA as correction function is varied for $p=3$, the region explored includes the peak CFL limit and the area where OOA is suspected to drop. What can be noted is that the analytically predicted decrease in order as $\iota_3 \rightarrow \infty$ is demonstrated numerically by $\mathrm{OOA}\rightarrow p$. Reassuringly, the order accuracy in the region of $\iota_3 \approx 10^{-3}$, where peak temporal stability is believed to lie, is shown to have $\mathrm{OOA} = p+1$ and hence the maximum CFL limit may be increased without an apparent loss in order accuracy. In particular, the GSFR correction function with peak temporal stability that was investigated in Fig.~\ref{fig:FRGS3_dispdiss} is found to recover the expected $p+1$ OOA. 
	
	The domain of $\mathbf{I}_p$ is sufficiently small for $p=3~\&~4$, therefore permitting an exhaustive search that couples numerical OOA calculation and analytic CFL calculation to find the maximum CFL limit. The results of this search are shown in Table~\ref{tab:FRGS_peak_CFL}, in all cases the OOA recovered was $p+1$, consistent with the expected order. It may be noted that in some cases a significant improvement upon the CFL limit may be made by comparison to Vincent~\etal~\cite{Vincent2011}. 
	
	\begin{figure}[tbhp]
		\centering
		\captionof{table}{Peak CFL of Order recovering GSFR}
		\label{tab:FRGS_peak_CFL}
		\begin{tabular}{|l|c|c|c|}
			\hline
			$p$ & RK Scheme & $\mathbf{I}_p$ & $\tau$ \\ \hline
			    & RK33 &$[1, 1.274\times10^{-3}, 1.438\times10^{-2}, 7.848\times10^{-3}]^T$ & 0.385 \\ 
			3   & RK44 &$[1, 2.069\times10^{-4}, 2.336\times10^{-3}, 2.336\times10^{-3}]^T$ & 0.390 \\
			    & RK55 &$[1, 6.952\times10^{-4},-6.158\times10^{-5}, 2.336\times10^{-3}]^T$ & 0.443 \\ \hline
			    & RK33 &$[1, 4.833\times10^{-4}, 2.336\times10^{-5}, -1.438\times10^{-4}, 2.637\times10^{-4}]^T$ & 0.431 \\
			4   & RK44 &$[1, 1.624\times10^{-3}, 2.637\times10^{-4}, -2.637\times10^{-4}, 2.637\times10^{-4}]^T$ & 0.430\\ 
			    & RK55 &$[1,1.624\times10^{-3},1.274\times10^{-5},-2.637\times10^{-4}, 8.859\times10^{-4}]^T$ & 0.354\\ \hline
		\end{tabular}
	\end{figure}
	
	\subsection{Linear Heterogeneous Advection}
	So far, the investigation has focused exclusively on the linear advection equation. However, for practical applications, non-linear conservation equations will be encountered. Other than the introduction of shock-waves, for which there treatment is still an important question, non-linear equations also introduce aliasing error due to the multiplication of polynomials. This results in a solution that lies beyond the spectral resolution of the grid. Further insight into the origin and effect of aliasing and order was investigated by Kravchenko and Moin~\cite{Kravchenko1997}. Previously, to investigate aliasing, a linear-heterogeneous flux function was used by Hesthaven and Warburton~\cite{Hesthaven2007} and Vincent~\etal~\cite{Vincent2010}. The flux function was defined as:
	\begin{equation}\label{eq:hest_flux}
		f(x;t) = \Big((1-x^2)^5 + 1\Big)u(x,t)
   	\end{equation}
   	However, when applied to a periodic domain $\mathbf{\Omega}=[-1,1]$ the flux function of Eq.~(\ref{eq:hest_flux}) is only $C^0$ continuous. Hence, it is proposed that it would be more suitable to use a flux function defined as:
   	\begin{equation}\label{eq:ali_flux_func}
   		f(x;t) = (\sin{(\pi x)}+2)u(x;t)
   	\end{equation}
   	The advantage of this flux function is that, although strictly linear in $u$, the spatial dependence of the flux function triggers aliasing errors. Therefore, this may give useful insight into GSFR's application to full non-linear problems, without the associated issues of shock formation. To understand the mechanism of production for aliasing error, consider the finite spatial Fourier series of an approximate solution $u(x,t)$:
   	\begin{equation}\label{eq:fourier_series}
   		u(x,t) = \sum^{N/2-1}_{n=-N/2} w_n(t)\exp{\Big(\frac{n\pi i x}{l}\Big)}
   	\end{equation}
   	where $w_n(t)$ are time dependent Fourier weights, and $l$ is the domain half-length. The half-length can be used to generalise Eq.~(\ref{eq:ali_flux_func}), such that, together with application of the product rule, one has:
   	\begin{equation}
   		\frac{\partial \big(\sin{(\frac{\pi x}{l}) + 2}\big)u}{\partial x} = \bigg(\sin{\Big(\frac{\pi x}{l}\Big)} + 2\bigg)\frac{\partial u}{\partial x} + \frac{\pi}{l}\cos{\Big(\frac{\pi x}{l}\Big)}u
	\end{equation}
	By subsequent application of the Fourier series of Eq.~(\ref{eq:fourier_series}):
	\begin{multline}
		\frac{\partial \big(\sin{(\frac{\pi x}{l}) + 2}\big)u}{\partial x} = \frac{\pi}{2l}\sum^{N/2-1}_{n=-N/2}  w_n(t) \Bigg( 4ni\exp{\Big(\frac{nx}{l}\Big)} \\
		  + \exp{\bigg(\frac{(n-1)\pi i x}{l}\bigg)}(1+ni) \\
		  + \exp{\bigg(\frac{(n+1)\pi i x}{l}\bigg)}(1-ni)\Bigg)
	\end{multline}
	Hence, aliasing is injected only by the highest frequency mode and will propagate down through the modes. This is behaviour is different to that expressed by flux function of Eq.~(\ref{eq:hest_flux}), where aliasing error is injected at all modes. 
	
	For numerical evaluation, a domain of $\mathbf{\Omega} = [-1,1]$ with periodic boundaries is used. The method of temporal integration used was low storage RK44, with $\tau$ set to be sufficiently small that the temporal integration had a negligible effect on the error. To illustrate the effect on the solution, the domain $L^2$ energy of the conserved variable is used, defined as:
   	\begin{equation}
		E(t) = \int_{\mathbf{\Omega}}u(x,t)^2dx
	\end{equation}   	 
	and the initial condition will be taken as:
	\begin{equation}
		u(x;t=0) = \sin{(4\pi x)}
	\end{equation}
	The effect of having a flux function that is $C^{\infty}$ continuous on $\mathbf{\Omega}$ is that the periodic boundaries are analogous to solving the same equation on an infinite domain. Hence, the time period of the solution can be found analytically to be $T=2/\sqrt{3}$. To evaluate the relative performance of correction functions, the $L^2$ energy error relative to some known value of the $L^2$ energy is compared for various correction functions. Due to the periodicity of the solution, the error at some time $nT,\:\: n\in\mathbb{N}$ may be straightforwardly calculated, \emph{i.e} analytical energy $E_a(nT) = 1$.
   
	\begin{figure}[tbhp]
		\centering
		\includegraphics[width=0.55\linewidth,trim= 0mm 0mm 0mm 0mm,clip=true]{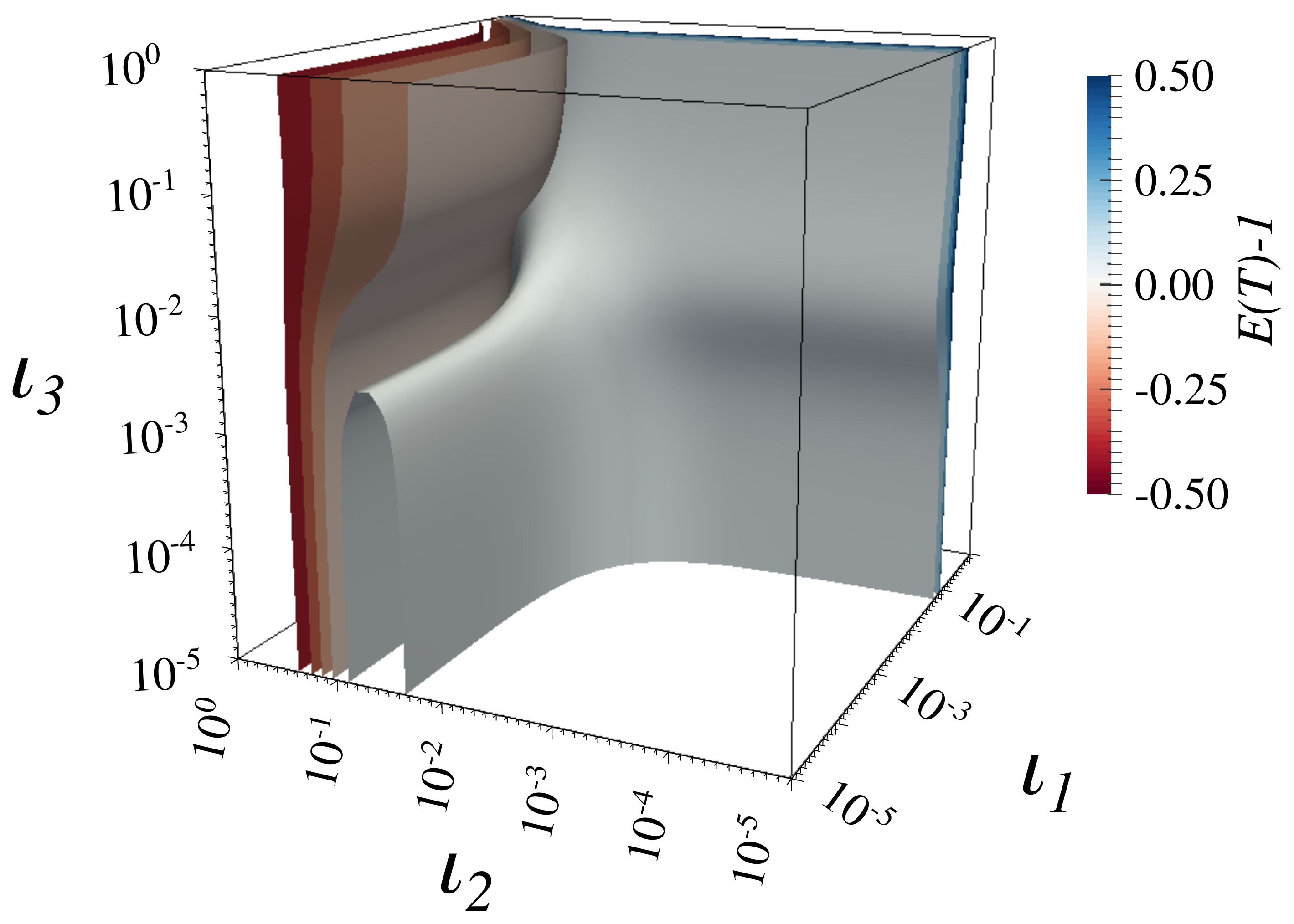}
		\caption{$L^2$ energy error for FR, $p=3$, using RK44 temporal integration with $\mathrm{CFL} = 0.06$, for $\alpha = 1$ (upwinded) $t=T=2/\sqrt{3}$.}
		\label{fig:FRGS3_aliL2}
	\end{figure}	
	 Figure~\ref{fig:FRGS3_aliL2} shows just such an error in the $L^2$ energy, for $p=3$ with upwinded interfaces. Central differenced interfaces are not shown as, for all correction functions, the scheme had gone unstable within $t=15T=30/\sqrt{3}$. This is due to the coupling of the  aliasing error to the zero dissipation associated with central difference. The result is that any error introduced into the solution does not become damped and hence instability arises. However, when considering the case of upwinded interfaces, the implicit dissipation in the scheme can be sufficient to stabilise the errors introduced through aliasing. It may also be noted that, from Fig.~\ref{fig:FRGS3_CFL}, the region of localised peak CFL number lies within a region of low aliasing error and hence correction functions in this region ($\iota_3 \approx1\times10^{-3},\iota_1 \approx 0$) may be good candidates for practical applications.
	
	\subsection{Euler's Equations}	
	The final investigation will be the application of GSFR  to the 3D Euler equations, defined as: 
	\begin{equation}
		\frac{\partial \mathbf{U}}{\partial t} + \nabla\cdot\mathbf{F} = 0
	\end{equation}
	where
	\begin{equation}
		\mathbf{F} = \begin{bmatrix}
		\mathbf{f}_1 & \mathbf{f}_2 & \mathbf{f}_3
		\end{bmatrix}^T
	\end{equation}
	\begin{equation}
		\mathbf{f}_1 = \begin{bmatrix}
		\rho u \\
		\rho u^2 + p \\
		\rho uv \\
		\rho uw \\
		u(\frac{1}{2}\rho\mathbf{v}\cdot\mathbf{v} + \rho e + p)
		\end{bmatrix}  \quad \mathbf{f}_2 = \begin{bmatrix}
		\rho v \\
		\rho uv \\
		\rho v^2 + p \\
		\rho vw \\
		v(\frac{1}{2}\rho\mathbf{v}\cdot\mathbf{v} + \rho e + p)
		\end{bmatrix}  \quad \mathbf{f}_3 = \begin{bmatrix}
		\rho w \\
		\rho uw \\
		\rho wv \\
		\rho w^2 + p \\
		w(\frac{1}{2}\rho\mathbf{v}\cdot\mathbf{v} + \rho e + p)
		\end{bmatrix}
	\end{equation}
	with $\mathbf{v} =[u,v,w]^T$. The case used to numerically investigate the performance of FR on Euler's Equations is the Isentropic Convecting Vortex (ICV) \cite{Cockburn1997}:
	\begin{align}
		u &= u_0 + \frac{\beta(y_c-y)}{2\pi}\exp{\bigg(\frac{1-r^2}{2}\bigg)} \\
		v &= v_0 + \frac{\beta(x-x_c)}{2\pi}\exp{\bigg(\frac{1-r^2}{2}\bigg)} \\
		w &= w_0 \\
		T & = 1 - \frac{(\gamma-1)\beta^2}{8\gamma\pi^2}\exp{\big(1-r^2\big)}
	\end{align}
	where, $\beta$ is the vortex strength, $x_c$ and $y_c$ are the grid centre coordinates, $\gamma$ is the ratio of specific heats, $r=\sqrt{x^2+y^2}$ and $T$ is temperature. The domain of the solution is taken as $\mathbf{\Omega}\in[-10,10]^3$ with periodic boundaries and the vortex strength $\beta=5$. For this case the common interface flux calculation used is  the Rusanov flux~\cite{Rusanov1961} with characteristic speed from Davis~\cite{Davis1988}.

	To extend the 1D scheme discussed here to three dimensions we will use a hexahedral grid. Consequently, a tensor product can be used to extend the spatial domain and the FR algorithm to incorporate other dimensions. This is the method outlined by Huynh~\cite{Huynh2007}.
	
	\begin{figure}[tbhp]
		\centering
		\includegraphics[width=0.9\linewidth]{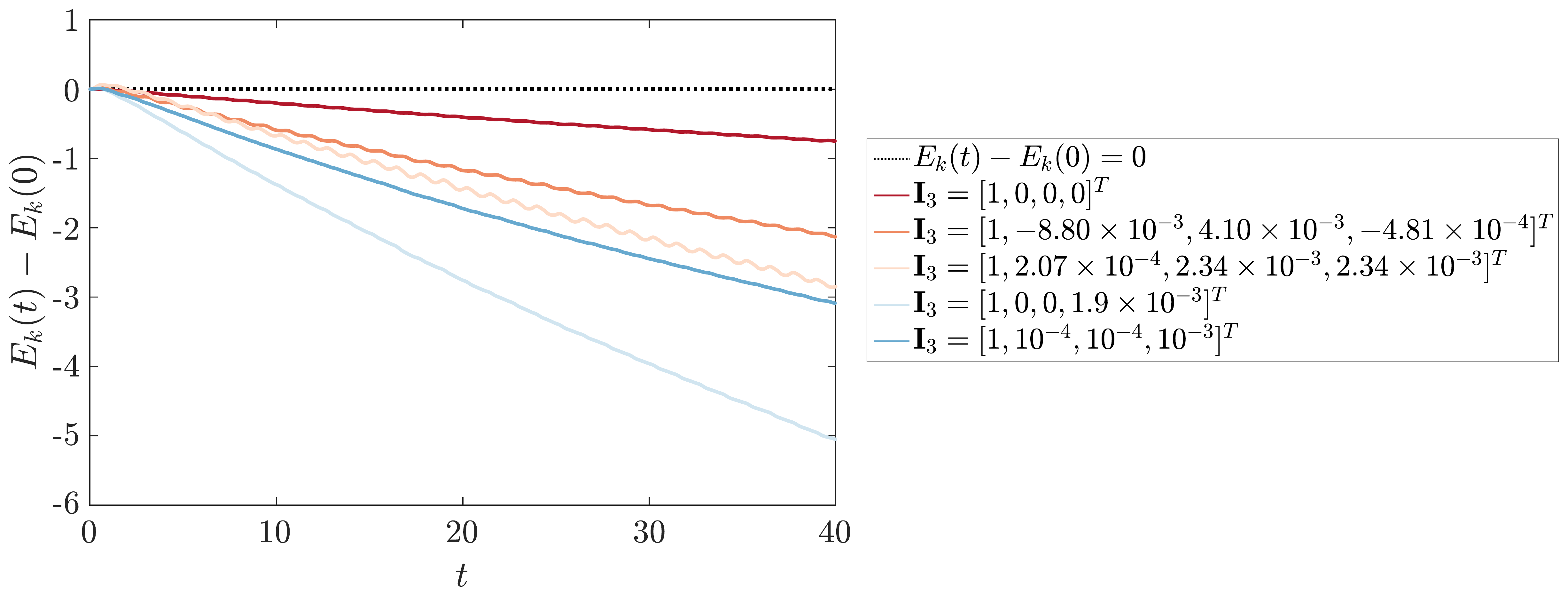}
		\caption{Error in the domain integrated kinetic energy of the ICV for FR.}
		\label{fig:FRGS3_ICV_EK}
	\end{figure}
	
	Figure~\ref{fig:FRGS3_ICV_EK} shows the error in the domain integrated kinetic energy of the ICV for $p=3$ using various correction functions. Here RK44 temporal integration with $\tau = 10^{-3}$ is used on $\mathbf{\Omega} = [-10,10]^3$ with $16\times16\times1$ elements. For the error calculation, $E_k(0)$ was calculated on a $256\times256\times1$ grid. Finally, the convective velocity was taken as $u_0=v_0=1, \:w_0 = 0$. 
	
	The correction functions applied here are the same as those whose dispersion and dissipation characteristics are shown in Fig.~\ref{fig:FRGS3_dispdiss}. All the correction functions seem to be stable when applied to the ICV. In all cases tested there is an small region, $0\leqslant t\leqslant 2$, where the kinetic energy increases above the initial value, this is likely due to the smoothing of the higher order discontinuities in the solution. The origin of these discontinuities is, although the domain is large, the solution is not strictly $C^{p}$ continuous on a periodic domain, and the effect of these discontinuities were seen to reduce as the domain is widened. This was investigated by Spiegel~\etal~\cite{Spiegel2015b}, where it was also shown that the periodic boundaries also caused instability through shear. Also visible in Fig.~\ref{fig:FRGS3_ICV_EK}, all the domain integrated error terms show a small degree of oscillation. This is due to the grid affecting integration accuracy and has a frequency that reflects this.
	
	Studying Fig.~\ref{fig:FRGS3_ICV_EK}, it can be clearly seen that DG gives the best performance. Although, it may have been expected that the GSFR CFL optimal correction function would give the best performance, as Fig.~\ref{fig:FRGS3_dispdiss} showed low dissipation at high wavenumbers. However, the grid used in this case is sufficient for the solution to be primarily formed of lower wavenumbers, where, as can be seen from the mid section of Fig.~\ref{fig:FRGS3_comp_I}, DG has the least dissipation.  An interesting comparison can be made between the performance of the optimal CFL correction functions for OSFR, ESFR, and GSFR. The OSFR case is clearly far more dissipative, from which it could be postulated that the lower order Sobolev terms present in ESFR and GSFR can lead to the cancelling out of unwanted dissipation.

\section{Conclusions}\label{sec:conclusions}
\begin{figure}
		\centering
		\includegraphics[width=0.4\linewidth]{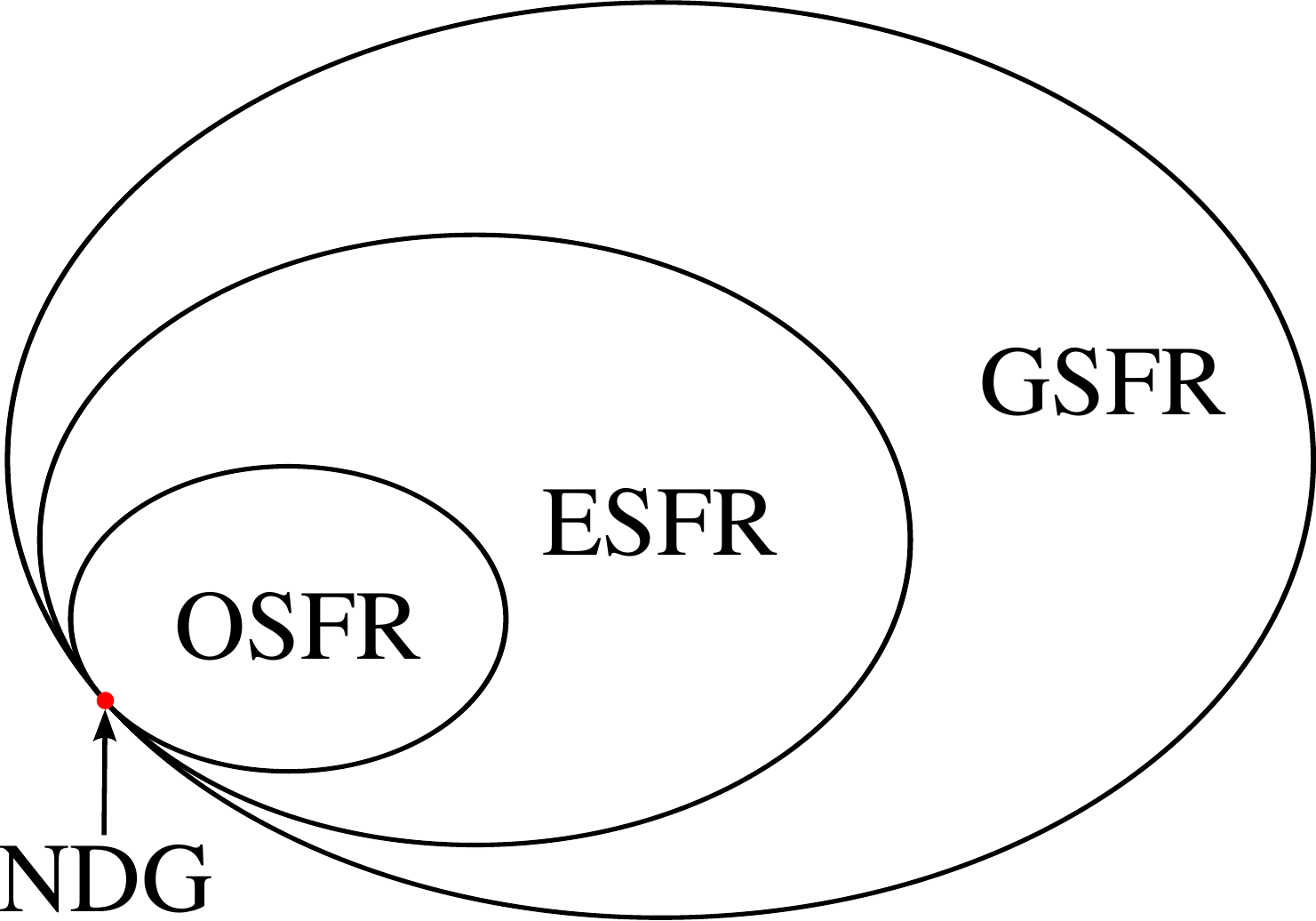}
		\caption{Diagram showing intersection of sets of correction functions.}
		\label{fig:FR_corr_vd}
	\end{figure}

	A new set of correction function has been presented, derived from the idea of energy stability in the broken Sobolev norm. This set of correction functions was shown to contain the OSFR and ESFR correction functions and it was further shown that unique functions could be obtained, that OSFR and ESFR were unable to produce. Hence, the intersection of the sets defining OSFR, ESFR and GSFR can be shown diagrammatically as in Fig.~\ref{fig:FR_corr_vd}.
	
	Study of the linear advection equation, via a von Neumann analysis, was then used to show that for $p=3$ a correction function could be found that led to an increase in the CFL number for RK44 temporal integration. It was also demonstrated that this correction function recovers numerically the expected OOA, with GSFR  also enabling the  recovery of functions from the lower order spaces. A one dimensional linear heterogeneous PDE was proposed to numerically investigate the performance of GSFR when applied to an equation that causes aliasing. The investigation showed that the region of high CFL limit for the case of $p=3$ was within a stable area of low error.  Finally, the solution of Euler's equations on the ICV test case was used to compare the relative performance when applied to fully non-linear PDEs. Of the correction functions tested all remained stable with reasonable performance, however it was found that DG gave the minimal decay in the kinetic energy for this test.

\section*{Acknowledgements}
\label{sec:ack}
	The support of the Engineering and Physical Sciences Research Council of the United Kingdom is gratefully acknowledged under the award reference 1750012. The author would like to acknowledge the useful discussions had with Rob Watson and his help in technical editing, language editing, and proofreading.

\section*{References}
\bibliographystyle{elsarticle-num}
\bibliography{library}


\clearpage
\begin{appendices}
\section{\textit{\textbf{L}\textsubscript{p}} Matrix Generation Algorithm}
\label{app:Lp}

	We will detail the basic algorithm for calculating the entries of the matrix $\mathbf{L}_p$. This is valid for the first p-2 rows of $\mathbf{L}_p$ and,  as can be seen form Eqs.(\ref{eq:corr_p2},\ref{eq:corr_p3},\ref{eq:corr_p4}), the penultimate two row enforce the boundary conditions of the left correction function.
	\begin{algorithm}
		\caption{Process for calculating the value of $\mathbf{L}_p[m-1][n]$ for a $\mathbf{I}_p=[\iota_0\dots\iota_p]$. The method of calculating $I_{i,m,n}$ is given in Algorithm~\ref{al:integral}.}
		\begin{algorithmic}
			\REQUIRE $p,m,n,\mathbf{I}_p$
			\STATE $L[m-1][n] \leftarrow 0$
			\FOR{$0\leqslant i \leqslant p$} 
				\STATE $L[m-1][n] \leftarrow L[m-1][n] + \iota_iI_{i,m,n}$
			\ENDFOR
			\FOR{$1\leqslant i \leqslant p$}
				\STATE $a \leftarrow \frac{(n+i)!(m+i)!}{2^{2i}(i!)^2(n-i)!(m-i)!}\big(1-(-1)^{n+m-2i}\big)$
				\STATE $L[m-1][n] \leftarrow L[m-1][n] - a\iota_i$ 
			\ENDFOR
			\RETURN $L[m-1][n]$
		\end{algorithmic}
	\end{algorithm}

	\begin{algorithm}
		\caption{Algorithm for calculating the integral of Eq.(\ref{eq:dmlegendre_int}) where $b_i(m,n)$ is defined by Eq.(\ref{eq:legendre_b}).}
		\begin{algorithmic}
			\REQUIRE $m,n,k$
			\STATE $I_{m,n,k} \leftarrow 0$
			\FOR{$0\leqslant i \leqslant \floor{\frac{n-m}{2}}$}
				\FOR{$0\leqslant j \leqslant \floor{\frac{k-m-1}{2}}$} 
					\STATE $c \leftarrow \frac{b_i(m,n)b_j(m+1,k)}{n+k-2(m+i+j)}\Big[1-(-1)^{n+k-2(m+i+j)}\Big]$
					\STATE $I_{m,n,k} \leftarrow I_{m,n,k} + c$
				\ENDFOR			
			\ENDFOR
			\RETURN $I_{m,n,k}$
		\end{algorithmic}
		\label{al:integral}
	\end{algorithm}
\section{Nomenclature}
\label{sec:nomenclauture}

	\begin{tabbing}
		  XXXXXX \= \kill
		  	\textit{Roman}\\
		  	$a_p$ \> $(2p)!/(2^p(p!)^2)$ \\
		  	$\mathbf{b}$ \> GSFR correction function boundary condition array\\
			$c(k)$ \> modified phase velocity at wavenumber $k$ \\
			$C^p$ \> $p^{\mathrm{th}}$ order differentiably continuous \\ 
			$\mathbf{C}_{+1}$ \> downwind cell FR matrix \\
			$\mathbf{C}_0$ \> centre cell FR matrix \\
			$\mathbf{C}_{-1}$ \> upwind cell FR matrix \\
			$\mathbf{D}$ \> first derivative matrix \\
			$E$ \> domain integrated energy of variable \\
			$E_k$ \> domain integrated kinetic energy \\
			$f$ \> flux variable in physical domain \\	
			$\mathbf{F}$ \> Array of Euler's equations flux vectors \\ 
			$h_l \:\mathrm{\&}\: h_r$ \> left and right correction functions\\
			$\mathbf{H}_p$ \> $p^{\mathrm{th}}$ order GSFR correction inverse matrix \\
			$g_l \:\mathrm{\&}\: g_r$ \> gradient of the left and right correction functions\\
			$J_i$ \> $i^{\mathrm{th}}$ cell Jacobian\\
			$k$ \> wavenumber \\
			$k_{nq}$ \> solution point Nyquist wavenumber, $(p+1)/\delta_j$\\
			$\hat{k}$ \> $k_{nq}$ normalised wavenumber, $[0,\pi]$ \\
			$\mathbf{K}$ \> ESFR correction matrix \\ 
			$l_i$ \> $i^{\mathrm{th}}$ Lagrange basis function \\
			$\mathbf{L}_p$ \> $p^{\mathrm{th}}$ order GSFR correction matrix \\
			$\mathbf{M}$ \> polynomial basis mass matrix \\ 
			$n_s$ \> number of discrete solution points in sub-domain\\ 
			$p$ \> solution polynomial order \\
			$\mathbf{Q}$ \> FR spatial discretisation operator matrix \\
			$\mathbf{R}$ \> FR spatial-temporal update matrix \\
			$u$ \> conserved variable in the physical domain \\
			$\mathbf{V}$ \> polynomial basis matrix \\
			$w_i$ \> $i^{\mathrm{th}}$ Fourier mode weight \\
			
			\\ \textit{Greek}\\
			$\alpha$ \> interface upwinding ratio ($\alpha = 1 \Rightarrow$ upwinded, $\alpha = 0.5 \Rightarrow$ central) \\
			$\delta_j$ \> mesh spacing, $x_j-x_{j-1}$\\
			$\epsilon_2$ \> domain averaged $L^2$ norm error \\
			$\iota$ \> OSFR correction function parameter \\
			$\iota_i$ \> $i^{\mathrm{th}}$ GSFR correction function parameter \\
			$\mathbf{I}_p$ \> $p^{\mathrm{th}}$ order array of GSFR correction function parameters, $\mathbf{I}_p = [\iota_0 \dots \iota_p]^T$\\
			$\kappa_i$ \> $i^{\mathrm{th}}$ ESFR correction parameter \\
			$\mathbf{\Lambda}$ \> $\mathrm{daig}(1,-1,1,-1\dots)$ \\
			$\xi$ \> transformed spatial variable \\
			$\rho(\mathbf{A})$ \> spectral radius of $\mathbf{A}$ \\ 
			$\tau$ \> time step \\
			$\psi_i$ \> $i^{\mathrm{th}}$ Legendre polynomial of the first kind \\ 
			$\mathbf{\Omega}$ \> solution domain \\
			$\mathbf{\Omega}_n$ \> $n^{\mathrm{th}}$ solution sub-domain \\
			$\hb{\Omega}$ \> reference sub-domain\\
			
			\\ \textit{Subscript}\\
			$\mathrm{\bullet}_l$ \> variable at left of cell\\
			$\mathrm{\bullet}_r$ \> variable at right of cell\\
			
			\\ \textit{Superscript}\\
			$\mathrm{\bullet}^T$ \> vector or matrix transpose\\
			$\mathrm{\bullet}^{\delta}$ \> discrete polynomial fit of value\\
			$\mathrm{\bullet}^{\delta C}$ \> correction to discontinuous discrete polynomial\\
			$\mathrm{\bullet}^{\delta D}$ \> discontinuous discrete polynomial fit of value\\
			$\mathrm{\bullet}^{\delta I}$ \> common value at interface based on discrete polynomial\\
			$\hat{\mathrm{\bullet}}$ \> variable transformed to reference domain\\
			$\tilde{\mathrm{\bullet}}$ \> variable transformed to Legendre basis\\
			
			\\ \textit{Operators}\\
			$\floor{x}$ \> floor function, giving largest integer less than or equal to $x$
		\end{tabbing}
\end{appendices}


\end{document}